\input ./style/arxiv-general.cfg
\documentclass[MSNbibl,nameyear,dvips]{arxstspdf}
\makeatletter
   \@ifpackageloaded{graphicx}{}{\usepackage{graphicx}}
\makeatother
\usepackage{flushend}
\usepackage{stfloats}

\volume{30}
\issue{3}
\pubyear{2015}
\firstpage{372}
\lastpage{390}
\doi{10.1214/15-STS520}
\docsubty{FLA}

\makeatletter

\newtheorem{theorem}{Theorem}
\newproclaim{definition}{Definition}
\newtheorem{corollary}{Corollary}
\newtheorem{lemma}{Lemma}
\newproclaim{remark}{Remark}
\newproclaim{assumption}{Assumption}
\makeatother

\begin{document}
\begin{frontmatter}

\title{Fourth Moments and Independent Component Analysis}%
\runtitle{Fourth moments and ICA}

\begin{aug}
\author[A]{\fnms{Jari}~\snm{Miettinen}\corref{}\ead[label=e1]{jari.p.miettinen@jyu.fi}},
\author[B]{\fnms{Sara}~\snm{Taskinen}\ead[label=e2]{sara.l.taskinen@jyu.fi}},
\author[C]{\fnms{Klaus}~\snm{Nordhausen}\ead[label=e3]{klaus.nordhausen@utu.fi}}
\and
\author[D]{\fnms{Hannu}~\snm{Oja}\ead[label=e4]{hannu.oja@utu.fi}}
\runauthor{Miettinen, Taskinen, Nordhausen and Oja}

\affiliation{University of Jyvaskyla, University of Jyvaskyla,
University of Turku and University of Turku}

\address[A]{Jari Miettinen is a~Postdoctoral Researcher, Department of Mathematics and
Statistics,
40014 University of Jyv\"askyl\"a, Jyv\"askyl\"a,
Finland \printead{e1}.}
\address[B]{Sara Taskinen is an Academy Research
Fellow, Department of Mathematics and
Statistics, 40014 University of Jyv\"askyl\"a, Jyv\"askyl\"a,
Finland \printead{e2}.}
\address[C]{Klaus Nordhausen is a~Senior Research Fellow, Department of Mathematics and
Statistics,
20014 University of Turku, Turku,   Finland \printead{e3}.}
\address[D]{Hannu Oja is a~Professor, Department of Mathematics and Statistics,
20014 University of Turku, Turku, Finland \printead{e4}.}
\end{aug}

%
\begin{abstract}
In independent component analysis it is assumed that the components of
the observed random
vector are linear combinations of latent independent random variables,
and the aim is then to find an estimate for a
transformation matrix back to these independent components. In the
engineering literature, there are
several traditional estimation procedures based on the use of fourth
moments, such as FOBI (fourth order blind identification), JADE (joint
approximate
diagonalization of eigenmatrices), and FastICA, but the statistical
properties of these estimates
are not well known. In this paper various independent component
functionals based on the
fourth moments are discussed in detail, starting with the corresponding
optimization problems, deriving the estimating equations and estimation
algorithms, and
finding asymptotic statistical properties of the estimates. Comparisons
of the asymptotic variances of the estimates in wide independent
component models
show that in most cases JADE and the symmetric version of FastICA
perform better than their competitors.
\end{abstract}

%
\end{frontmatter}

\section{Introduction}\label{sec1}

In his system of frequency curves, \citet{Pearson1895} identified
different types of distributions,
and the classification was based on the use of the standardized third
and fourth moments.
A measure of degree of kurtosis for the distribution of $x$ was defined as
\[
\beta=\frac{E ([x-E(x)]^4 )} { [ E ([x-E(x)]^2
) ]^2}\quad\mbox{or}\quad\kappa=\beta-3,
\]
and \citet{Pearson1905} called the distribution platykurtic,
leptokurtic, or mesokurtic depending on
the value of $\kappa$. In the case of the normal distribution ($\kappa
=0$, mesokurtic) Pearson also
considered the probable error of $\hat\kappa$. Later, kurtosis was
generally understood simply as a
property which is measured by $\kappa$, which has raised questions such
as ``Is kurtosis really
peakedness?''; see, for example, \citet{Darlington1970}.
\citet{vanZwet1964} proposed kurtosis orderings for
symmetrical distributions, and \citet{Oja1981} defined measures of
kurtosis as functionals which (i)
are invariant under linear transformations and (ii) preserve the van
Zwet partial ordering. Most of
the measures of kurtosis, including $\beta$, can be written as a ratio
of two scale measures.
Recently, robust measures of kurtosis also have been proposed and
considered in the literature;
see, for example, \citet{Brys2006}.

It is well known that the variance of the sample mean depends on the
population variance only,
but the variance of the sample variance depends also on the shape of
the distribution through
$\beta$. The measure $\beta$ has been used as an indicator of the
bimodality, for example, in identifying
clusters in the data set (\citeauthor{Pena2001}, \citeyear{Pena2001})
or as a general indicator
for non-Gaussianity, for example, in
testing for normality or in independent component analysis (\citeauthor
{Hyvarinen1999}, \citeyear{Hyvarinen1999}). Classical tests
for the normality are based on the standardized third and fourth
moments. See also \citet{DeCarlo1997}
for the meaning and use of kurtosis.

The concept and measures of kurtosis have been extended to the
multivariate case as well. The
classical skewness and kurtosis measures by \citet{Mardia1970}, for
example, combine in a natural
way the third and fourth moments of a standardized multivariate
variable. Mardia's measures are
invariant under affine transformations, that is, the $p$-variate random
variables $\mathbf{x}$ and
$\mathbf{A}\mathbf{x}+\mathbf{b}$ have the same skewness and
kurtosis values for all
full-rank $p\times p$ matrices
$\mathbf{A}$ and for all $p$-vectors $\mathbf{b}$. For similar
combinations of
the standardized third and fourth
moments, see also \citet{Morietal1993}. Let next $\mathbf{V}_1$
and $\mathbf{V}_2$ be two $p\times p$ affine
equivariant scatter matrices (functionals); see \citet{Huber1981} and
\citet{Maronna1976} for
early contributions on scatter matrices. Then, in the invariant
coordinate selection (ICS) in
\citet{Tyler2009}, one finds an affine transformation matrix
$\mathbf{W}$
such that
\[
\mathbf{W} \mathbf{V}_1\mathbf{W}'=
\mathbf{I}_p \quad\mbox {and}\quad \mathbf{W}\mathbf{V}_2
\mathbf{W}'=\mathbf{D},
\]
where $\mathbf{D}$ is a diagonal matrix with diagonal elements in decreasing
order. The transformed
$p$ variables are then presented in a new invariant coordinate system,
and the diagonal elements
in $\mathbf{D}$, that is, the eigenvalues of $\mathbf
{V}_1^{-1}\mathbf{V}_2$, provide
measures of multivariate
kurtosis. This procedure is also sometimes called the generalized
principal component analysis and
has been used to find structures in the data. See
\citet{Caussinus1993}, \citet{Critchley2007},
\citet{Ilmonenetal2010a}, \citet{Pena2010}, and \citet
{Nordhausenetal2011b}. For the tests for
multinormality based on these ideas, see \citet
{Kankainenetal2007}. In
independent component
analysis, certain fourth moment matrices are used together with the
covariance matrix in a similar
way to find the transformations to independent components [FOBI
by \citet{Cardoso1989} and JADE
by \citet{CardosoSouloumiac1993}]. See also \citet{Oja2006}.

In this paper, we consider the use of univariate and multivariate
fourth moments in independent
component analysis (ICA). The basic independent component (IC) model
assumes that the observed
components of $\mathbf{x}_i=(x_{i1},\ldots,x_{ip})'$ are linear combinations
of latent independent
components of $\mathbf{z}_i=(z_{i1},\ldots,z_{ip})'$. Hence, the
model can be
written as
\[
\mathbf{x}_i=\bolds{\mu}+\bolds{\Omega}\mathbf{z}_i,
\quad i=1, \ldots,n,
\]
where the full rank $p \times p$ matrix $\bolds{\Omega}$ is called the
mixing matrix and
$\mathbf{z}_1,\ldots,\mathbf{z}_n$ is a random sample from a
distribution with
independent components such
that $E(\mathbf{z}_i)=\mathbf{0}$ and $\operatorname{Cov}(\mathbf
{z}_i)=\mathbf{I}_p$.
Similarly to the
model of elliptically symmetric
distributions, the IC model is a semiparametric model, as the marginal
distributions of the
components of $\mathbf{z}$ are left fully unspecified except for the first
two moments. For the
identifiability of the parameters, one further assumes that at most one
of the components has a
normal distribution. Notice also that $\bolds{\Omega}$ and $\mathbf
{z}$ are
still confounded in the
sense that the order and signs of the components of $\mathbf{z}$ are not
uniquely defined.
The location center, the $p$-vector $\bolds{\mu}$, is usually
considered a
nuisance parameter,
since the main goal in independent component analysis is, based on a
$p\times n$ data matrix
$\mathbf{X}=(\mathbf{x}_1,\dots,\mathbf{x}_n)$, to find an estimate
for an unmixing
matrix $\mathbf{W}$ such that
$\mathbf{W}\mathbf{x}$ has independent components. Note that all unmixing
matrices $\mathbf{W}$ can be written
as $\mathbf{C}\bolds{\Omega}^{-1}$, where each row and each column
of the
$p\times p$ matrix $\mathbf{C}$ has
exactly one nonzero element.

The population quantity to be estimated is first defined as an
independent component functional
$\mathbf{W}(F)$. The estimate $\mathbf{W}(F_n)$, also denoted by
$\mathbf{W}(\mathbf{X})$,
is then obtained by applying the functional
to the empirical distribution $F_n$ of $\mathbf{X}=(\mathbf
{x}_1,\ldots,\mathbf{x}_n)$.
In the engineering literature, several
estimation procedures based on the fourth moments, such as FOBI (fourth
order blind
identification) (\citeauthor{Cardoso1989}, \citeyear{Cardoso1989}),
JADE (joint approximate
diagonalization of eigenmatrices)
(\citeauthor{CardosoSouloumiac1993}, \citeyear
{CardosoSouloumiac1993}), and FastICA (\citeauthor{Hyvarinen1999},
\citeyear{Hyvarinen1999}),
have been proposed and widely used. In these
approaches the marginal distributions are separated using various
fourth moments. On the other hand, the estimators
by \citet{ChenBickel2006} and \citet{Samworth2012} only
need the
existence of the first moments and
rely on efficient nonparametric estimates of the marginal densities.
Efficient estimation methods based
on residual signed ranks and residual ranks have been developed
recently by \citet{IlmonenPaindaveine2011} and
\citet{Hallin2015}. For a parametric model with a marginal Pearson
system approach, see
\citet{KarvanenKoivunen2002}.

This paper describes in detail the independent component functionals
based on fourth moments through
corresponding optimization problems and estimating equations, provides
fixed-point algorithms and the limiting
statistical properties of the estimates, and specifies the needed
assumptions. Also, a
wide comparison study of the estimates is carried out. As far as we
know, most of the
results in the paper are new, including the asymptotical properties of
the JADE estimate. The
asymptotical properties of the FOBI estimate have been derived earlier
in \citet{Ilmonenetal2010a}.
The limiting variances and the limiting multinormality of the
deflation-based version of the FastICA
estimate have been studied in \citet{Ollila2010}
and \citet{Nordhausenetal2011}, respectively.

\section{Notation and Preliminary Results}

Throughout the paper, we use the following notation. First write, for
independent
$z_{ik}$, $k=1,\dots,p$,
\begin{eqnarray*}
E(z_{ik}) &=& 0,\quad E \bigl(z_{ik}^2 \bigr)=1,
\quad E \bigl(z_{ik}^3 \bigr)=\gamma_k \quad
\mbox{and}
\\
E \bigl(z_{ik}^4 \bigr) &=& \beta_k,
\end{eqnarray*}
and
\[
\kappa_k=\beta_k-3,\quad\pi_k=
\operatorname{sign}( \kappa_k)\quad\mbox{and}\quad\operatorname{Var}
\bigl(z_{ik}^3 \bigr)=\sigma_k^2.
\]
As seen later, the limiting distributions of the unmixing matrix
estimates based on fourth moments
depend on the joint limiting distribution of
%
%
\begin{eqnarray}
\sqrt{n} \hat s_{kl} &=& n^{-1/2} \sum
_{i=1}^n z_{ik}z_{il},
\nonumber
\\[-8pt]
\label{rkl}
\\[-8pt]
\nonumber
\sqrt{n} \hat r_{kl} & =& n^{-1/2} \sum
_{i=1}^n \bigl(z_{ik}^3-
\gamma_k \bigr)z_{il}
\end{eqnarray}
and
\[
\sqrt{n} \hat r_{mkl}= n^{-1/2} \sum
_{i=1}^n z_{im}^2z_{ik}z_{il},
\]
for distinct $k,l,m=1,\ldots,p$. If the eighth moments of $z_i$ exist,
then the joint limiting
distribution of $\sqrt{n} \hat s_{kl}$, $\sqrt{n} \hat r_{kl}$, and
$\sqrt{n} \hat r_{mkl}$ is a multivariate normal distribution with
marginal zero means. The
nonzero variances and covariances are
\begin{eqnarray*}
\operatorname{Var}(\sqrt{n} \hat s_{kl}) &=& 1,\quad\operatorname
{Var}(\sqrt{n} \hat r_{kl})= \sigma_k^2,
\\
\operatorname{Var}(\sqrt{n} \hat r_{mkl}) &=& \beta_m,
\end{eqnarray*}
and
\begin{eqnarray*}
\operatorname{Cov}(\sqrt{n} \hat s_{kl},\sqrt{n} \hat
r_{kl})&=& \beta_k,
\\
\operatorname{Cov}(\sqrt{n} \hat r_{kl},\sqrt{n} \hat
r_{lk})&=&\beta_k\beta_l,
\end{eqnarray*}
and
\begin{eqnarray*}
\operatorname{Cov}(\sqrt{n} \hat s_{kl},\sqrt{n} \hat
r_{mkl})&=&1,
\\
\operatorname{Cov}( \sqrt{n} \hat r_{kl},\sqrt{n} \hat
r_{mkl})&=& \beta_k \quad\mbox{and}
\\
\operatorname{Cov}(\sqrt{n} \hat r_{lk},\sqrt{n} \hat
r_{mkl}) &=& \beta_l.
\end{eqnarray*}

We also often refer to the following sets of $p\times p$ transformation
matrices:
\begin{enumerate}[5.]
\item[1.] $\mathcal{D}=\{\operatorname{diag}(d_1,\ldots,d_p):
d_1,\ldots,d_p>0 \}$
(heterogeneous rescaling),
\item[2.] $\mathcal{J}=\{\operatorname{diag}(j_1,\ldots,j_p):
j_1,\ldots,j_p=\pm1 \}$
(heterogeneous sign changes),
\item[3.] $\mathcal{P}=\{\mathbf{P}: \mathbf{P} \mbox{ is a
permutation matrix}
\}$,
\item[4.] $\mathcal{U}=\{\mathbf{U}: \mathbf{U} \mbox{ is an orthogonal
matrix}\}$,
\item[5.] $\mathcal{C}=\{\mathbf{C}: \mathbf{C}=\mathbf{P}\mathbf
{J}\mathbf{D}, \mathbf{P}\in
\mathcal{P}, \mathbf{J}\in\mathcal{J},
\mathbf{D}\in\mathcal{D} \}$.
\end{enumerate}
Next, let $\mathbf{e}_i$ denote a $p$-vector with $i$th element one and
other elements zero, and define
$\mathbf{E}^{ij}=\mathbf{e}_i\mathbf{e}_j'$, $i,j=1,\dots,p$, and
\begin{eqnarray*}
\mathbf{J}_{p,p} &=& \sum_{i=1}^p
\sum_{j=1}^p \mathbf{E}^{ij}
\otimes\mathbf{E}^{ij} =\operatorname{vec} (\mathbf{I}_p)
\operatorname{vec} (\mathbf{I}_p)',
\\
\mathbf{K}_{p,p} &=& \sum_{i=1}^p
\sum_{j=1}^p \mathbf{E}^{ij}
\otimes\mathbf{E}^{ji},
\\
\mathbf{I}_{p,p} &=& \sum_{i=1}^p
\sum_{j=1}^p \mathbf{E}^{ii}
\otimes\mathbf{E}^{jj} =\mathbf {I}_{p^2}\quad\mbox{and}
\\
\mathbf{D}_{p,p} &=& \sum_{i=1}^p
\mathbf{E}^{ii} \otimes\mathbf{E}^{ii}.
\end{eqnarray*}
Then, for any $p\times p$ matrix $\mathbf{A}$,
$\mathbf{J}_{p,p} \operatorname{vec}(\mathbf{A}) = \operatorname
{tr}(\mathbf{A})\cdot\break \operatorname{vec}(\mathbf{I}_p)$,
$\mathbf{K}_{p,p} \operatorname{vec}(\mathbf{A}) = \operatorname
{vec}(\mathbf{A}')$, and
$\mathbf{D}_{p,p} \operatorname{vec}(\mathbf{A}) = \operatorname
{vec}(\operatorname{diag}
(\mathbf{A}))$.
The matrix $\mathbf{K}_{p,p}$ is sometimes called a commutation
matrix. For
a symmetric nonnegative
definite matrix $\mathbf{S}$, the matrix $\mathbf{S}^{-1/2}$ is taken
to be
symmetric\vspace*{1.3pt} and nonnegative
definite and to satisfy $\mathbf{S}^{-1/2} \mathbf{S} \mathbf
{S}^{-1/2}=\mathbf{I}_p$.

\section{Independent Component Model and~Functionals}

\subsection{Independent Component (IC) Model}

Throughout the paper, our $p$-variate observations $\mathbf{x}_1,\dots
,\mathbf{x}_n$ follow the independent component
(IC) model
%
%
\begin{eqnarray}
\label{ICmodel} \mathbf{x}_i&= & \bolds{\mu}+\bolds{\Omega}
\mathbf{z}_i,\quad i=1,\ldots,n,
\end{eqnarray}
where $\bolds{\mu}$ is a mean vector, $\bolds{\Omega}$ is a
full-rank $p\times
p$ mixing matrix, and
$\mathbf{z}_1,\ldots,\mathbf{z}_n$ are independent and identically
distributed
random vectors from a
$p$-variate distribution such that:
%

\begin{assumption}\label{a1}
The components $z_{i1},\ldots,z_{ip}$ of $\mathbf{z}_i$ are independent.
\end{assumption}

%
\begin{assumption}\label{a2}
Second moments exist, $E(\mathbf{z}_i)=\mathbf{0}$ and $E(\mathbf
{z}_i\mathbf{z}_i')=\mathbf{I}_p$.
\end{assumption}

%
\begin{assumption}\label{a3}
At most one of the components $z_{i1},\ldots,z_{ip}$ of $\mathbf
{z}_i$ has a
normal distribution.
\end{assumption}

If the model is defined using Assumption~\ref{a1} only, then the mixing
matrix $\bolds{\Omega}$ is not
well-defined and can at best be identified only up to the order, the
signs, and heterogenous
multiplications of its columns. Assumption~\ref{a2} states that the
second moments exist, and
$E(\mathbf{z}_i)=\mathbf{0}$ and $E(\mathbf{z}_i\mathbf
{z}_i')=\mathbf{I}_p$ serve as
identification constraints for $\bolds{\mu}$
and the scales of the columns of
$\bolds{\Omega}$. Assumption~\ref{a3} is needed, as, for example, if
$\mathbf{z}\sim N_2(\mathbf{0},\mathbf{I}_2)$,
then also $\mathbf{U}\mathbf{z}\sim N_2(\mathbf{0},\mathbf{I}_2)$
for all orthogonal $\mathbf{U}$ and the independent
components are not well-defined.
Still, after these three assumptions, the order and signs of the
columns of $\bolds{\Omega}$ remain
unidentified, but one can identify the set of the standardized
independent components
$\{\pm z_{i1},\dots,\pm z_{ip}\}$, which is naturally sufficient for
practical data analysis.

One of the key results in independent component analysis is the following.

%
\begin{theorem}
\label{thm1}
Let $\mathbf{x}=\bolds{\mu}+\bolds{\Omega}\mathbf{z}$ be an
observation from
an IC model
with mean vector $\bolds{\mu}$
and covariance matrix $\bolds{\Sigma}=\bolds{\Omega}\bolds{\Omega
}'$, and write
$\mathbf{x}_{st}=\bolds{\Sigma}^{-1/2}(\mathbf{x}-\bolds{\mu})$
for the standardized random variable. Then $\mathbf{z}=\mathbf
{U}\mathbf{x}_{st}$ for
some orthogonal matrix
$\mathbf{U}=(\mathbf{u}_1,\ldots,\mathbf{u}_p)'$.
\end{theorem}

The result says that, starting with standardized observations $\mathbf
{x}_{st}$, one only has to search for an unknown
$\mathbf{U}\in\mathcal{U}$ such that $\mathbf{U}\mathbf{x}_{st}$
has independent
components. Thus,
after estimating $\bolds{\Sigma}$, the estimation problem can be
reduced to
the estimation problem
of an orthogonal matrix $\mathbf{U}$ only.

\subsection{Independent Component (IC) Functionals}

Write next $\mathbf{X}=(\mathbf{x}_1,\ldots,\mathbf{x}_n)$ for a
random sample from
the IC model~(\ref{ICmodel}) with the
cumulative distribution function (c.d.f.) $F_{\mathbf{x}}$. As mentioned
in the \hyperref[sec1]{Introduction}, the aim of independent component analysis (ICA) is
to find an estimate
of some unmixing matrix $\mathbf{W}$ such that $\mathbf{W} \mathbf
{x}_i$ has
independent components.
It is easy to see that all unmixing matrices can be written as $\mathbf
{W}=\mathbf{C}\bolds{\Omega}^{-1}$ for
some $\mathbf{C}\in\mathcal{C}$. The population quantity, which we
wish to
estimate, is defined as
the value of an independent component functional $\mathbf{W}(F)$ at the
distribution of $F_{\mathbf{x}}$.

%
\begin{definition}\label{ICfunc}
The $p\times p$ matrix-valued functional $\mathbf{W}(F)$ is said to be an
independent component (IC)
functional if (i) $\mathbf{W}(F_{\mathbf{x}}) \mathbf{x}$ has
independent components
in the IC
model~(\ref{ICmodel}) and
(ii) $\mathbf{W}(F_{\mathbf{x}})$ is affine equivariant in the sense that
\begin{eqnarray*}
&& \bigl\{ \bigl(\mathbf{W}(F_{\mathbf{A}\mathbf{x}+\mathbf
{b}})\mathbf{A}\mathbf{x}
\bigr)_1,\ldots , \bigl( \mathbf{W}(F_{\mathbf{A}\mathbf{x}+\mathbf{b}})\mathbf{A}
\mathbf {x} \bigr)_p \bigr\}
\\
&&\quad= \bigl\{ \pm \bigl(\mathbf{W}(F_{\mathbf{x}})\mathbf {x}
\bigr)_1,\dots,\pm \bigl(\mathbf{W}(F_{\mathbf{x}})\mathbf {x}
\bigr)_p \bigr\}
\end{eqnarray*}
for all nonsingular $p \times p$ matrices $\mathbf{A}$ and for all
$p$-vectors $\mathbf{b}$.
\end{definition}

Notice that in the independent component model, $\mathbf{W}(F_{\mathbf
{x}})\mathbf{x}$
does not depend on the
specific choices of $\mathbf{z}$ and $\bolds{\Omega}$, up to the
signs and the
order of the components. Notice also that, in the condition (ii), any
c.d.f. $F$ is allowed to be used as an argument of $\mathbf{W}(F)$. The
corresponding sample version
$\mathbf{W}(F_n)$ is then obtained when the IC functional is applied
to the
empirical distribution function $F_n$ of $\mathbf{X}=(\mathbf
{x}_1,\ldots,\mathbf{x}_n)$. We also sometimes write $\mathbf
{W}(\mathbf{X})$ for the sample version.
Naturally, the estimator is then also affine
equivariant in the sense that, for all nonsingular $p \times p$
matrices $\mathbf{A}$ and for all $p$-vectors $\mathbf{b}$, $\mathbf
{W}(\mathbf{A}\mathbf{X}+\mathbf{b}\mathbf{1}_n')\mathbf{A}\mathbf
{X}=\mathbf{P}\mathbf{J}\mathbf{W}(\mathbf{X})\mathbf{X}$ for some
$\mathbf{J}\in
\mathcal{J}$ and $\mathbf{P}\in\mathcal{P}$.

%
\begin{remark}
As mentioned before, if $\mathbf{W}$ is an unmixing matrix, then so is
$\mathbf{C}\mathbf{W}$ for all $\mathbf{C}\in\mathcal{C}$,
and we then have a whole set of matrices $\{ \mathbf{C}\mathbf{W}:
\mathbf{C}\in
\mathcal{C}\}$ equivalent to $\mathbf{W}$.
To find a unique representative in the class, it is often required that
$\operatorname{Cov}(\mathbf{C}\mathbf{W}\mathbf{x})=\mathbf{I}_p$ but
still the order and signs of the rows remain unidentified. Of course,
the assumption
on the existence of second moments may sometimes be thought to be too
restrictive. For alternative ways to
identify the unmixing matrix, see then \citet{ChenBickel2006},
\citet{IlmonenPaindaveine2011}, and
\citet{Hallin2015}, for example. For a general discussion on this
identification problem, see also
\citet{ErikssonKoivunen2004}.
\end{remark}

\section{Univariate Kurtosis and Independent Component Analysis}

\subsection{Classical Measures of Univariate Skewness and Kurtosis}

Let first $x$ be a univariate random variable with mean value $\mu$ and
variance $\sigma^2$.
The standardized variable is then $z=(x-\mu)/\sigma$, and classical
skewness and kurtosis measures are
the standardized third and fourth moments, $\gamma=E(z^3)$ and $\beta
=E(z^4)$. For symmetrical distributions,
$\gamma=0$, and for the normal distribution, $\kappa=\beta-3=0$.
For a random sample $x_1,\dots,x_n$ from a univariate distribution, write
\begin{eqnarray*}
\mu_j &=& E \bigl((x_i-\mu)^j \bigr)\quad
\mbox{and}
\\
m_j &=& n^{-1} \sum_{i=1}^n
(x_i-\bar x)^j,\quad j=2,3,4.
\end{eqnarray*}
Then the limiting distribution of
$
\sqrt{n} (m_2-\mu_2,m_3-\mu_3,m_4-\mu_4)'
$
is a 3-variate normal distribution with mean vector zero and covariance
matrix with the $(i,j)$ element
\begin{eqnarray*}
&& \mu_{i+j+2} - \mu_{i+1}\mu_{j+1}-(i+1)
\mu_i\mu_{j+2}
\\
&&\quad{}-(j+1)\mu_{i+2}\mu_j+ (i+1) (j+1)
\mu_i\mu_j\mu_2,
\end{eqnarray*}
$i,j=1,2,3$. See Theorem~2.2.3.B in \citet{Serfling1980}.
Then in the symmetric case with $\mu_2=1$, for example,
\begin{eqnarray*}
&&\sqrt{n} \pmatrix{ m_2-1
\cr
\vspace*{0pt} m_3
\cr
\vspace*{0pt} m_4-\mu_4 }
\\
&&\quad\to_d N_3 \left( \pmatrix{ 0
\cr
0
\cr
0},\right.
\\
&&\left.{}\qquad{}\pmatrix{ \mu_4-1 & 0 & \mu_6-
\mu_4
\cr
\vspace*{0pt} 0 & \mu_6-6\mu_4+9 & 0
\cr
\vspace*{0pt} \mu_6-\mu_4 & 0 & \mu_8-
\mu_4^2 } \right).
\end{eqnarray*}
If the observations come from $N(0,1)$, we further obtain
\[
\sqrt{n} \pmatrix{ m_2-1
\cr
\vspace*{0pt} m_3
\cr
\vspace*{0pt} m_4-3} \to_d N_3 \left( \pmatrix{ 0
\cr
\vspace*{0pt} 0
\cr
\vspace*{0pt} 0}, \pmatrix{ 2 & 0 & 12
\cr
\vspace*{0pt} 0 & 6
& 0
\cr
\vspace*{0pt} 12 & 0 & 96} \right).
\]
The classical skewness and kurtosis statistics, the natural estimates
of $\gamma$ and $\beta$, are
$\hat\gamma= {m_3}/{m_2^{3/2}}$ and $\hat\beta= {m_4}/{m_2^2}$,
and then
\begin{eqnarray*}
\sqrt{n} \hat\gamma&=& \sqrt{n} m_3+o_P(1)\quad
\mbox{and}
\\
\sqrt{n} \hat\kappa&=& \sqrt{n} (\hat\beta-3)
\\
&=&\sqrt{n} (m_4-3)-6 \sqrt{n} (m_2-1)+o_P(1)
\end{eqnarray*}
and we obtain, in the general $N(\mu,\sigma^2)$ case, that
\begin{eqnarray*}
\sqrt{n}  \pmatrix{ \hat\gamma
\cr
\vspace*{0pt} \hat\kappa} &=& \sqrt{n}
\pmatrix{ \hat\gamma
\cr
\vspace*{0pt} \hat\beta-3 }
\\
& \to_d & N_2 \left( \pmatrix{ 0
\cr
0}, \pmatrix{ 6 & 0
\cr
0 & 24 } \right).
\end{eqnarray*}
Consider next $p$-variate observations coming from an IC model. The
important role of the fourth moments
is stated in the following:
%

\begin{theorem}\label{key2}
Let the components of $\mathbf{z}=(z_1,\ldots, \break z_p)'$ be
independent and
standardized so that $E(\mathbf{z})=\mathbf{0}$ and $\operatorname
{Cov}(\mathbf{z})=\mathbf{I}_p$, and
assume that
at most one of the kurtosis values $\kappa_i=E(z_i^4)-3$, $i=1,\ldots
,p$, is zero. Then the following inequalities hold true:
\begin{longlist}
\item[(i)]
\begin{eqnarray*}
&& \bigl|E \bigl( \bigl(\mathbf{u}' \mathbf{z} \bigr)^4
\bigr)-3 \bigr|
\\
&&\quad\leq\max \bigl\{ \bigl|E \bigl(z_1^4 \bigr)-3 \bigr|,
\ldots, \bigl|E \bigl(z_p^4 \bigr)-3 \bigr| \bigr\}
\end{eqnarray*}
for all $\mathbf{u}$ such that $\mathbf{u}' \mathbf{u}=1$. The
equality holds only if
$\mathbf{u}=\mathbf{e}_i$ for $i$
such that $|E(z_i^4)-3|= \max\{|E(z_1^4)-3|,\ldots, |E(z_p^4)-3|
\}$,
and
\item[(ii)]
\begin{eqnarray*}
&& \bigl|E \bigl[ \bigl(\mathbf{u}_1' \mathbf{z}
\bigr)^4 \bigr]-3 \bigr|+\cdots+ \bigl|E \bigl[ \bigl(\mathbf{u}_p'
\mathbf{z} \bigr)^4 \bigr]-3 \bigr|
\\
&&\quad\le \bigl|E \bigl[z_1^4 \bigr]-3\bigr|+ \cdots+\bigr|E
\bigl[z_p^4 \bigr]-3 \bigr|
\end{eqnarray*}
for all orthogonal matrices $\mathbf{U}=(\mathbf{u}_1,\dots,\mathbf
{u}_p)'$. The
equality holds only if $\mathbf{U}=\mathbf{J}\mathbf{P}$ for some
$\mathbf{J}\in\mathcal
{J}$ and $\mathbf{P}\in\mathcal{P}$.
\end{longlist}
\end{theorem}

For the first part of the theorem, see Lemma~2 in \citet{BugrienKent2005}.
The theorem suggests natural strategies and algorithms in search for
independent components.
It was seen in Theorem~\ref{thm1} that in the IC model $\mathbf
{x}_{st}=\mathbf{U}\mathbf{z}$
with an orthogonal $\mathbf{U}=(\mathbf{u}_1,\ldots,\mathbf{u}_p)$.
The first part of
Theorem~\ref{key2} then shows
how the components can be found one by one just by repeatedly maximizing
\[
\bigl|E \bigl( \bigl(\mathbf{u}_k' \mathbf{x}_{st}
\bigr)^4 \bigr)-3 \bigr|,\quad k=1,\ldots,p
\]
(projection pursuit approach), and the second part of Theorem~\ref
{key2} implies that the same components may be
found simultaneously by maximizing
\[
\bigl|E \bigl[ \bigl(\mathbf{u}_1' \mathbf{x}_{st}
\bigr)^4 \bigr]-3 \bigr|+\cdots+ \bigl|E \bigl[ \bigl(\mathbf{u}_p'
\mathbf{x}_{st} \bigr)^4 \bigr]-3 \bigr|.
\]
In the engineering literature, these two approaches are well known and
important special cases of the so-called
deflation-based FastICA and symmetric FastICA; see, for example,
\citet{HyvarinenKarhunenOja2001}. The statistical
properties of these two estimation procedures will now be considered in detail.

\subsection{Projection Pursuit Approach---Deflation-Based FastICA}

Assume that $\mathbf{x}$ is an observation from an IC model~(\ref{ICmodel})
and let again
$\mathbf{x}_{st}=\bolds{\Sigma}^{-1/2}(\mathbf{x}-\bolds{\mu})$
be the standardized
random variable.
Theorem~\ref{key2}(i) then suggests the following projection pursuit approach
in searching for the independent components.

%
\begin{definition}
The deflation-based projection pursuit (or deflation-based
FastICA) functional
is$\mathbf{W}(F_{\mathbf{x}})=\mathbf{U}\bolds{\Sigma
}^{-1/2}$, where
$\bolds{\Sigma}=\operatorname{Cov}(\mathbf{x})$
and the rows of an orthogonal matrix $\mathbf{U}=(\mathbf{u}_1,\dots
,\mathbf{u}_p)'$
are found one by one
by maximizing
\[
\bigl|E \bigl( \bigl(\mathbf{u}_k' \mathbf{x}_{st}
\bigr)^4 \bigr)-3 \bigr|
\]
under the constraint that $\mathbf{u}_k' \mathbf{u}_k=1$ and $\mathbf
{u}_j' \mathbf{u}_k=0$, $j=1,\ldots,k-1$.
\end{definition}

It is straightforward to see that $\mathbf{W}(F_{\mathbf{x}})$ is affine
equivariant. In the independent
component model~(\ref{ICmodel}), $\mathbf{W}(F_{\mathbf{x}})\mathbf
{x}$ has
independent components if
Assumption~\ref{a3} is replaced by the following stronger assumption.
%

\begin{assumption}\label{a4}
The fourth moments of $\mathbf{z}$ exist, and at most one of the kurtosis
values $\kappa_k$, $k=1,\ldots,p$, is zero.
\end{assumption}

Thus, under this assumption, $\mathbf{W}(F)$ is an independent component
(IC) functional. Based on Theorem~\ref{key2}(i), the functional
then finds the independent components in such an order that
\[
\bigl|E \bigl( \bigl(\mathbf{u}_1' \mathbf{x}_{st}
\bigr)^4 \bigr)-3 \bigr|\ge\cdots\ge \bigl|E \bigl( \bigl(\mathbf{u}_p'
\mathbf{x}_{st} \bigr)^4 \bigr)-3 \bigr|.
\]
The solution order is unique if the kurtosis values are distinct.

The Lagrange multiplier technique can be used to obtain the estimating
equations for
$\mathbf{U}=(\mathbf{u}_1,\ldots,\mathbf{u}_p)'$.
This is done in \citet{Ollila2010} and \citet
{Nordhausenetal2011} and
the procedure is the following.
After finding $\mathbf{u}_1,\dots,\mathbf{u}_{k-1}$, the solution
$\mathbf{u}_k$ thus
optimizes the Lagrangian function
\[
L(\mathbf{u}_k ,\bolds{\theta}_k)= \bigl|E \bigl( \bigl(
\mathbf {u}_k' \mathbf{x}_{st}
\bigr)^4 \bigr)-3 \bigr|- \sum_{j=1}^k
\theta_{kj} \bigl(\mathbf{u}_j'
\mathbf{u}_k-\delta _{jk} \bigr),
\]
where $\bolds{\theta}_k=(\theta_{k1},\dots,\theta_{kk})'$ is the
vector of
Lagrangian multipliers and
$\delta_{jk}=1$ $(0)$ as $j=k$ ($j\ne k$) is the Kronecker delta.
Write
\[
\mathbf{T}(\mathbf{u})=E \bigl[ \bigl(\mathbf{u}'
\mathbf{x}_{st} \bigr)^3 \mathbf{x}_{st} \bigr].
\]
The solution for $\mathbf{u}_k$ is then given by the $p+k$ equations
\begin{eqnarray*}
4 \pi_k \mathbf{T}(\mathbf{u}_k)- \sum
_{j=1}^{k-1} \theta_{kj}
\mathbf{u}_j -2\theta_{kk} \mathbf{u}_k &=&
\mathbf{0}\quad\mbox{and}
\\
\mathbf{u}_j' \mathbf{u}_k &=&
\delta_{jk},\quad j=1, \ldots,k,
\end{eqnarray*}
where $\pi_k=\operatorname{sign}(\kappa_k)$.
One then first finds the solutions for the Lagrange coefficients in
$\theta_k$, and substituting
these results into the first $p$ equations, the following result is obtained.
%

\begin{theorem}\label{est-eq-fastICA}
Write $\mathbf{x}_{st}=\bolds{\Sigma}^{-1/2}(\mathbf{x}-\bolds{\mu
})$ for the
standardized random vector, and
$\mathbf{T}(\mathbf{u})=E [(\mathbf{u}' \mathbf{x}_{st} )^3\cdot\break \mathbf{x}_{st} ]$. The
orthogonal matrix
$\mathbf{U}=(\mathbf{u}_1,\dots,\mathbf{u}_p)'$ solves the
estimating equations
\begin{eqnarray}
&& \bigl( \mathbf{u}_k' \mathbf{T}(
\mathbf{u}_k) \bigr) \mathbf {u}_k= \Biggl(
\mathbf{I}_p-\sum_{j=1}^{k-1}
\mathbf{u}_j\mathbf{u}_j' \Biggr) \mathbf{T}(
\mathbf {u}_k),
\nonumber
\\
\eqntext{\displaystyle k=1,\ldots,p.}
\end{eqnarray}
\end{theorem}

The theorem suggests the following fixed-point algorithm for the
deflation-based solution.
After finding $\mathbf{u}_1,\dots,\mathbf{u}_{k-1}$, the following
two steps are
repeated until convergence
to get $\mathbf{u}_k$:
\begin{eqnarray*}
&&\mbox{Step 1: } \quad\mathbf{u}_{k} \leftarrow
\Biggl(I_p-\sum_{j=1}^{k-1}
\mathbf{u}_j\mathbf{u}_j' \Biggr) \mathbf{T}(
\mathbf{u}_{k}),
\\
&&\mbox{Step 2:}\quad\mathbf{u}_{k} \leftarrow\Vert
\mathbf{u}_{k} \Vert^{-1} \mathbf{u}_{k}.
\end{eqnarray*}

The deflation-based estimate $\mathbf{W}(\mathbf{X})$ is obtained as
above but by
replacing the population
quantities by the
corresponding empirical ones. Without loss of generality, assume next that
$|\kappa_1|\ge\cdots\ge|\kappa_p|$.
First note that, due to the affine equivariance of the estimate,
$\mathbf{W}(\mathbf{X})=\mathbf{W}(\mathbf{Z})\bolds{\Omega
}^{-1}$. In the efficiency studies,
it is therefore sufficient
to consider $\hat{\mathbf{W}}=\mathbf{W}(\mathbf{Z})$
and the limiting distribution of $\sqrt{n} (\hat{\mathbf{W}}-\mathbf
{I}_p)$ for
a sequence $\hat{\mathbf{W}}$ converging in probability to $\mathbf{I}_p$.
As the empirical and population criterion functions
\begin{eqnarray*}
D_n(\mathbf{u}) &=& \Biggl\vert n^{-1} \sum
_{i=1}^n \bigl(\mathbf{u}'
\mathbf{x}_{st,i} \bigr)^4-3 \Biggr\vert \quad\mbox{and}
\\
D(\mathbf{u}) &=& \bigl\vert E \bigl[ \bigl(\mathbf{u}' \mathbf{z}
\bigr)^4 \bigr]-3 \bigr\vert
\end{eqnarray*}
are continuous and $\sup_{\mathbf{u}' \mathbf{u}=1} |D_n(\mathbf
{u})-D(\mathbf{u})|\to_P
0$, one can choose a sequence of solutions such that $\hat{\mathbf
{u}}_1\to
_P \mathbf{e}_1$ and similarly for $\hat{\mathbf{u}}_2,\ldots,\hat
{\mathbf{u}}_{p-1}$.
Further, then also $\hat{\mathbf{W}}=\hat{\mathbf{U}}\hat{\mathbf
{S}}^{-1/2}\to_P
\mathbf{I}_p$. One can
next show that the limiting distribution of $\sqrt{n} (\hat{\mathbf
{W}}-\mathbf{I}_p)$ is obtained if we
only know the joint limiting distribution of $\sqrt{n} (\hat{\mathbf
{S}}-\mathbf{I}_p)$ and $\sqrt{n} \operatorname{off}(\hat{\mathbf{R}})$,
where $\hat{\mathbf{S}}=(\hat s_{kl})$ is the sample covariance
matrix, $\hat{\mathbf{R}}=(\hat r_{kl})$ is given in (\ref{rkl}),
and $\operatorname{off}(\hat{\mathbf{R}})=\hat{\mathbf
{R}}-\operatorname{diag}(\hat{\mathbf{R}})$. We then have
the following results; see also \citet{Ollila2010}, \citet
{Nordhausenetal2011}.

%
\begin{theorem}
\label{asymp6}
Let $\mathbf{Z}=(\mathbf{z}_1,\ldots,\mathbf{z}_n)$ be a random
sample from a
distribution with finite eighth
moments and satisfying the Assumptions \ref{a1}, \ref{a2}, and \ref{a4}
with $|\kappa_1|\ge\cdots\ge|\kappa_p|$.
Then there exists a sequence of solutions such that $\hat{\mathbf
{W}}\to_{P}
\mathbf{I}_p$ and
\begin{eqnarray}
\sqrt{n} \hat{w}_{kl}&=& -\sqrt{n} \hat{w}_{lk}-\sqrt{n}
\hat{s}_{kl}+o_{P}(1), \quad l<k,
\nonumber
\\
\sqrt{n} (\hat{w}_{kk}-1) &=& -1/2\sqrt{n} (\hat{s}_{kk}-1)+o_{P}(1)
\quad\mbox{and}
\nonumber
\\
\sqrt{n} \hat{w}_{kl}&=& \frac{\sqrt{n} \hat r_{kl}-(\kappa
_k+3)\sqrt
{n} \hat{s}_{kl}} {\kappa_k}+o_{P}(1),
\nonumber
\\
\eqntext{\displaystyle l>k.}
\end{eqnarray}
\end{theorem}

%
\begin{corollary}
Under the assumptions of Theorem~\ref{asymp6}, the limiting
distribution of $\sqrt{n} \operatorname{vec} (\hat{\mathbf
{W}}-\mathbf{I}_p)$ is
a~multivariate normal with zero mean vector and componentwise variances
\begin{eqnarray*}
\operatorname{ASV}(\hat{w}_{kl})&=& \frac{\sigma_l^2 -(\kappa
_l+3)^2 } {\kappa
_l^2}+1,\quad
\kappa_l \ne0, l<k,
\\
\operatorname{ASV}(\hat{w}_{kk}) &=& (\kappa_k+2)/4 \quad
\mbox{and}
\\
\operatorname{ASV}(\hat{w}_{kl})&=& \frac{\sigma_k^2 -(\kappa
_k+3)^2 } {\kappa_k^2 },\quad
\kappa_k\ne0, l>k.
\end{eqnarray*}
\end{corollary}

%
\begin{remark}
Projection pursuit is used to reveal structures in the
original data
by selecting interesting low-dimensional orthogonal projections of
interest. This is done, as above,
by maximizing the value of an objective function (projection index).
The term ``projection pursuit''
was first launched by \citet{FriedmanTukey1974}. \citet{Huber1985}
considered projection indices
with heuristic arguments that a projection is the more interesting, the
less normal it is. All his
indices were ratios of two scale functionals, that is, kurtosis
functionals, with the classical kurtosis
measure as a special case. He also discussed the idea of a recursive
approach to find subspaces.
\citet{Pena2001} used the projection pursuit algorithm with the
classical kurtosis index for finding
directions for cluster identification. For more discussion on the
projection pursuit approach, see also \citet{JonesSibson1987}.
\end{remark}

%
\begin{remark}
In the engineering literature,\break \citet{HyvarinenOja1997}
were the first
to propose the procedure based
on the fourth moments, and later considered an extension with a choice
among several alternative projection
indices (measures of non-Gaussianity). The approach is called
deflation-based or one-unit FastICA and it is
perhaps the most popular approach for the ICA problem in engineering
applications. Note that the estimating
equations in Theorem~\ref{est-eq-fastICA} and the resulting fixed-point
algorithm do not fix the
order of the components (the order is fixed by the original definition)
and, as seen in
Theorem~\ref{asymp6}, the limiting distribution of the estimate depends
on the order in which the
components are found. Using this property, \citet{Nordhausenetal2011}
proposed a two-stage version
of the deflation-based FastICA method with a chosen projection index
that finds the components in an optimal efficiency order.
Moreover, \citet{Miettinenetal2014a} introduced an adaptive two-stage
algorithm that (i) allows one to use different projection indices
for different components and (ii) optimizes the order in which the
components are extracted.
\end{remark}

\subsection{Symmetric Approach---Symmetric FastICA}

In the symmetric approach, the rows of the matrix $\mathbf{U}$ are found
simultaneously,
and we have the following:

%
\begin{definition}
The symmetric projection pursuit (or symmetric\vspace*{1pt} fastICA)
functional is
$\mathbf{W}(F_{\mathbf{x}})=\mathbf{U}\bolds{\Sigma}^{-1/2}$,
where $\bolds{\Sigma}=\operatorname{Cov}(\mathbf{x})$ and
$\mathbf{U}=(\mathbf{u}_1,\dots,\mathbf{u}_p)'$ maximizes
\[
\bigl|E \bigl( \bigl(\mathbf{u}_1' \mathbf{x}_{st}
\bigr)^4 \bigr)-3 \bigr|+ \cdots+ \bigl|E \bigl( \bigl(\mathbf{u}_p'
\mathbf{x}_{st} \bigr)^4 \bigr)-3 \bigr|
\]
under the constraint that $\mathbf{U} \mathbf{U}'=\mathbf{I}_p$.
\end{definition}

This optimization procedure is called symmetric FastICA in the signal
processing community.
The functional $\mathbf{W}(F_{\mathbf{x}})$ is again affine
equivariant. Based on
Theorem~\ref{key2}(ii),
in the IC model with Assumption~\ref{a4} the maximizer is unique up to
the order
and signs of the rows of $\mathbf{U}$, that is,
\[
\{ z_{1},\ldots,z_{p}\} = \bigl\{ \pm\mathbf{u}_1'
\mathbf {x}_{st},\ldots, \pm\mathbf{u}_p'
\mathbf{x}_{st} \bigr\}.
\]

As in the deflation-based case, we use the Lagrange multiplier
technique to obtain the
matrix $\mathbf{U}$. The Lagrangian function to be optimized is now
\begin{eqnarray*}
L(\mathbf{U},\bolds{\Theta})&=&\sum_{k=1}^p
\bigl|E \bigl( \bigl(\mathbf{u}_k' \mathbf{x}_{st}
\bigr)^4 \bigr)-3 \bigr|- \sum_{k=1}^p
\theta_{kk} \bigl(\mathbf{u}'_k
\mathbf{u}_k-1 \bigr)
\\
&&{}-\sum_{j=1}^{p-1}\sum
_{k=j+1}^p\theta_{jk}\mathbf{u}'_j
\mathbf{u}_k,
\end{eqnarray*}
where the symmetric matrix $\bolds{\Theta}=(\theta_{jk})$ contains all
$p(p+1)/2$ Lagrangian multipliers.\vspace*{1.5pt}
Write again
$\mathbf{T}(\mathbf{u})=E((\mathbf{u}'\mathbf{x}_{st})^3 \mathbf{x}_{st})$.
Then the solution $\mathbf{U}=(\mathbf{u}_1,\ldots,\mathbf
{u}_p)'$ satisfies
\begin{eqnarray}
&& 4 \pi_k T(\mathbf{u}_k)=2\theta_{kk}
\mathbf{u}_k+\sum_{j<k}
\theta_{jk}\mathbf{u}_j+\sum_{j>k}
\theta_{kj}\mathbf {u}_j,
\nonumber
\\
\eqntext{\displaystyle k=1,\ldots,p,}
\end{eqnarray}
and
\[
\mathbf{U}\mathbf{U}'=\mathbf{I}_p.
\]
Solving $\theta_{jk}$ and using the fact that $\theta_{jk}=\theta_{kj}$
give $\pi_k \mathbf{u}_j'
\mathbf{T}(\mathbf{u}_k)=\pi_j \mathbf{u}_k' \mathbf{T}(\mathbf
{u}_j)$, $j,k=1,\dots,p$,
and\vspace*{1.8pt} we get the following estimating equations.
%

\begin{theorem}
Let $\mathbf{x}_{st}=\bolds{\Sigma}^{-1/2}(\mathbf{x}-\bolds{\mu
})$ be the
standardized
random vector from the IC model~(\ref{ICmodel}),
$\mathbf{T}(\mathbf{u})=E((\mathbf{u}' \mathbf{x}_{st})^3\mathbf{x}_{st})$,
$
\mathbf{T}(\mathbf{U})=(\mathbf{T}(\mathbf{u}_1),\dots,\mathbf
{T}(\mathbf{u}_p))'$ and $\bolds{\Pi
}=\operatorname{diag}
(\pi_1,\dots,\pi_p)$.
The estimating equations for the symmetric solution $\mathbf{U}$ are
\[
\mathbf{U} \mathbf{T}(\mathbf{U})' \bolds{\Pi}=\bolds{\Pi }
\mathbf{T}(\mathbf{U})\mathbf{U}' \quad\mbox{and}\quad\mathbf {U}
\mathbf{U}' =\mathbf{I}_p.
\]
\end{theorem}

For the computation of $\mathbf{U}$, the above estimating equations suggest
a fixed-point algorithm with
the updating step
\[
\mathbf{U}\leftarrow\bolds{\Pi}\mathbf{T} \bigl(\mathbf{T}'
\mathbf{T} \bigr)^{-1/2}.
\]

The symmetric version estimate $\mathbf{W}(\mathbf{X})$ is obtained by
replacing the population quantities by their corresponding
empirical\vspace*{1pt}
ones in the estimating equations.
Write again $\hat{\mathbf{W}}=\mathbf{W}(\mathbf{Z})$ and let $\hat
{\mathbf{S}}=(\hat s_{kl})$
and $\hat{\mathbf{R}} =(\hat r_{kl})$ be as in (\ref{rkl}). Then we
have the
following:


\begin{theorem}
\label{asymp7}
Let $\mathbf{Z}=(\mathbf{z}_1,\ldots, \mathbf{z}_n)$ be a random
sample from a
distribution of $\mathbf{z}$ satisfying
the Assumptions \ref{a1}, \ref{a2}, and \ref{a4} with bounded eighth
moments. Then there is a
sequence of solutions such that $\hat{\mathbf{W}}\to_{P}\mathbf
{I}_p$ and
\begin{eqnarray*}
&&\!\!\sqrt{n} (\hat{w}_{kk}-1)
\\
&&\!\!\!\quad= - \frac{1}2 \sqrt{n} (\hat {s}_{kk}-1)+o_{P}(1)
\quad\mbox{and}
\\
&&\!\!\sqrt{n} \hat{w}_{kl}
\\
&&\!\!\!\quad= \frac{\sqrt{n} \hat r_{kl}\pi_k- \sqrt{n}
\hat r_{lk}\pi_l-(\kappa_k\pi_k+3\pi_k-3\pi_l)
\sqrt{n} \hat{s}_{kl} } {|\kappa_k| + |\kappa_l| }
\\
&&\qquad{}+o_{P}(1), \quad k\neq l,
\end{eqnarray*}
where $\pi_k=\operatorname{sign}(\kappa_k)$.
\end{theorem}

%
\begin{corollary}
\label{sficalim}
Under the assumptions of Theorem~\ref{asymp7}, the limiting
distribution of
$\sqrt{n} \operatorname{vec} (\hat{\mathbf{W}}-\mathbf{I}_p)$ is a
multivariate normal
with zero mean vector and
componentwise variances
\begin{eqnarray}
\operatorname{ASV}(\hat{w}_{kk}) &=& (\kappa_k+2)/4\quad
\mbox{and}
\nonumber
\\
\operatorname{ASV}(\hat{w}_{kl}) &=& \frac{\sigma_k^2+\sigma
_l^2-\kappa_k^2-6(\kappa
_k+\kappa_l)-18 }{
(|\kappa_k|+|\kappa_l| )^2},
\nonumber
\\
\eqntext{\displaystyle k\neq l.}
\end{eqnarray}
\end{corollary}

%
\begin{remark}
The symmetric FastICA approach with other choices of projection indices
was proposed in the engineering
literature by \citet{Hyvarinen1999}. The computation of symmetric
FastICA estimate was done, as
in our approach, by running $p$ parallel one-unit algorithms, which
were followed by a matrix
orthogonalization step. A generalized symmetric FastICA algorithm that
uses different projection
indices for different components was proposed by~\citet
{KoldovskyTichavskyOja2006}. The asymptotical
variances of generalized symmetric FastICA estimates were derived
in~\citet{TichavskyKoldovskyOja2006} under the assumption of symmetric
independent component
distributions.
\end{remark}

\section{Multivariate Kurtosis and Independent Component Analysis}

\subsection{Measures of Multivariate Skewness and~Kurtosis}

Let $\mathbf{x}$ be a $p$-variate random variable with mean vector
$\bolds
{\mu}$
and covariance matrix
$\bolds{\Sigma}$, and $\mathbf{x}_{st}= \bolds{\Sigma}^{-1/2}
(\mathbf{x}-\bolds{\mu})$. All
the standardized third and fourth moments can now be collected into
$p\times p^2$ and
$p^2\times p^2$ matrices
\begin{eqnarray*}
\bolds{\gamma} &=& E \bigl(\mathbf{x}_{st}'\otimes \bigl(
\mathbf {x}_{st} \mathbf{x}_{st}' \bigr) \bigr)
\quad\mbox{and}
\\
\bolds{\beta} &=& E \bigl( \bigl(\mathbf{x}_{st} \mathbf{x}_{st}'
\bigr)\otimes \bigl(\mathbf{x}_{st} \mathbf{x}_{st}'
\bigr) \bigr).
\end{eqnarray*}
Unfortunately, these multivariate measures of skewness and kurtosis are
not invariant under affine
transformations: The transformation $\mathbf{x}\to\mathbf{A}\mathbf
{x}+\mathbf{b}$
induces, for some unspecified
orthogonal matrix $\mathbf{U}$, the transformations
\begin{eqnarray*}
\mathbf{x}_{st} &\to& \mathbf{U} \mathbf{x}_{st},\quad
\bolds{ \gamma}\to\mathbf{U} \bolds{\gamma} \bigl(\mathbf{U}'\otimes
\mathbf{U}' \bigr)\quad\mbox{and}
\\
\bolds{\beta} &\to& (\mathbf{U}\otimes\mathbf{U}) \bolds{\beta } \bigl(
\mathbf{U}'\otimes\mathbf{U}' \bigr).
\end{eqnarray*}

Notice next that, for any $p\times p$ matrix $\mathbf{A}$,
%
%
\begin{eqnarray}
\mathbf{G}(\mathbf{A}) &=& E \bigl(\mathbf{x}_{st}
\mathbf{x}_{st}' \mathbf{A} \mathbf{x}_{st}
\bigr) \quad\mbox{and}
\nonumber
\\[-8pt]
\label{34moments}
\\[-8pt]
\nonumber
\mathbf{B}(\mathbf{A}) &=& E \bigl( \mathbf{x}_{st}
\mathbf{x}_{st}' \mathbf{A}\mathbf{x}_{st}
\mathbf{x}_{st}' \bigr)
\end{eqnarray}
provide selected $p$ and $p^2$ linear combinations of the third and
fourth moments as
$\operatorname{vec}(\mathbf{G}(\mathbf{A}))=\gamma\operatorname
{vec}(\mathbf{A})$
and $\mbox
{vec}(\mathbf{B}(\mathbf{A}))=\beta\operatorname{vec}(\mathbf{A})$.
Further, the elements of matrices
\[
\mathbf{G}^{ij}= \mathbf{G} \bigl(\mathbf{E}^{ij} \bigr)
\quad\mbox {and}\quad\mathbf{B}^{ij}=\mathbf{B} \bigl(
\mathbf{E}^{ij} \bigr),\quad i,j=1,\dots,p,
\]
list all possible third and fourth moments. Also,
\[
\mathbf{G}=\mathbf{G}(\mathbf{I}_p)=\sum
_{i=1}^p \mathbf{G}^{ii}\quad\mbox{and}
\quad\mathbf{B}=\mathbf{B}(\mathbf {I}_p)=\sum
_{i=1}^p \mathbf{B}^{ii}
\]
appear to be natural measures of multivariate skewness and kurtosis. In
the independent component
model we then have the following straightforward result.


\begin{theorem}
\label{B}
At the distribution of $\mathbf{z}$ with independent components,
$E(\mathbf{z})=\mathbf{0}$, $\operatorname{Cov}(\mathbf{z})=\mathbf{I}_p$,
and $\kappa_i=E(z_i^4)-3$, $i=1,\ldots,p$:
\begin{eqnarray}
\bolds{\beta}&=&\sum_{i=1}^p
\kappa_i \bigl(\mathbf{E}^{ii}\otimes\mathbf{E}^{ij}
\bigr)+ \mathbf{I}_{p,p}+\mathbf{J}_{p,p}+\mathbf{K}_{p,p},
\nonumber
\\
\mathbf{B}^{ij}&=&\sum_{k=1}^p
\kappa_k \bigl(\mathbf{E}^{kk}\mathbf{E}^{ij}
\mathbf{E}^{kk} \bigr)+\mathbf{E}^{ij}+\mathbf{E}^{ji}+
\operatorname{tr} \bigl(\mathbf {E}^{ij} \bigr)\mathbf{I}_p,
\nonumber
\\
\eqntext{\displaystyle i,j=1,\ldots,p\quad\mbox{and}}
\\
\mathbf{B} &=& \sum_{i=1}^p (
\kappa_i+p+2)\mathbf{E}^{ii}.
\nonumber
\end{eqnarray}
\end{theorem}

%
\begin{remark}
The standardized third and fourth moments have been used as building
bricks for invariant
multivariate measures of skewness and kurtosis. The classical skewness
and kurtosis measures
by~\citet{Mardia1970} are
\[
E \bigl( \bigl(\mathbf{x}_{st}' \tilde{
\mathbf{x}}_{st} \bigr)^3 \bigr) \quad\mbox{and}\quad
\operatorname{tr}(\mathbf{B} )=E \bigl( \bigl(\mathbf{x}_{st}'
\mathbf{x}_{st} \bigr)^2 \bigr),
\]
whereas \citet{Morietal1993} proposed
\begin{eqnarray}
\Vert\mathbf{G}\Vert^2 &=& E \bigl(\mathbf{x}_{st}'
\mathbf {x}_{st} \mathbf{x}_{st}' \tilde{
\mathbf{x}}_{st} \tilde{\mathbf{x}}_{st}' \tilde{
\mathbf{x}}_{st} \bigr) \quad \mbox{and}
\nonumber
\\
\operatorname{tr}(\mathbf{B} ) &=& E \bigl( \bigl(\mathbf{x}_{st}'
\mathbf{x}_{st} \bigr)^2 \bigr),
\nonumber
\end{eqnarray}
where $\mathbf{x}_{st}$ and $\tilde{\mathbf{x}}_{st}$ are
independent copies of
$\mathbf{x}_{st}$\break (\citeauthor{Morietal1993}, \citeyear
{Morietal1993}). (The invariance follows as $\mathbf{x}\to\mathbf{A}
\mathbf{x} +\mathbf{b}$ induces
$\mathbf{x}_{st}\to\mathbf{U} \mathbf{x}_{st}$ for some orthogonal
$\mathbf{U}$.) The
sample statistics can then be
used to test multivariate normality, for example. For their limiting
distributions under the normality
assumption, see, for example, \citet{Kankainenetal2007}. For other
extensions of multivariate skewness
and kurtosis and their connections to skewness and kurtosis measures
above, see \citet{Kollo2008}
and \citet{KolloSrivastava2005}.
In Sections~\ref{sFOBI} and~\ref{sJADE}, we first use $\mathbf{B}$ alone
and then all $\mathbf{B}^{ij}$,
$i,j=1,\ldots,p$, together to find solutions to the independent
component problem. In the signal
processing literature, these approaches are called FOBI (fourth order
blind identification) and
JADE (joint approximate diagonalization of eigenmatrices), correspondingly.
\end{remark}

\subsection{Use of Kurtosis Matrix $\mathbf{B}$---FOBI}
\label{sFOBI}

The independent component functional based on the covariance matrix
$\bolds{\Sigma}$ and the kurtosis
matrix $\mathbf{B}$ defined in~(\ref{34moments}) is known as FOBI (fourth
order blind identification)
(\citeauthor{Cardoso1989}, \citeyear{Cardoso1989}) in the engineering
literature. It is one of the
earliest approaches to the
independent component problem and is defined as follows.
%

\begin{definition}
The FOBI functional is $\mathbf{W}(F_{\mathbf{x}})=\mathbf{U}\bolds
{\Sigma}^{-1/2}$,
where $\bolds{\Sigma}=\operatorname{Cov}(\mathbf{x})$ and the rows of
$\mathbf{U}$ are the eigenvectors of $\mathbf{B}=E (\mathbf
{x}_{st}\mathbf{x}_{st}'
\mathbf{x}_{st} \mathbf{x}_{st}')$.
\end{definition}

First recall that, in the independent component model, $\mathbf
{x}_{st}=\mathbf{U}'\mathbf{z}$ for some orthogonal
$\mathbf{U}$. This implies that
\[
\mathbf{B}=E \bigl(\mathbf{x}_{st}\mathbf{x}_{st}'
\mathbf{x}_{st} \mathbf{x}_{st}' \bigr)=
\mathbf{U}' E \bigl(\mathbf{z} \mathbf{z}' \mathbf{z}
\mathbf{z}' \bigr) \mathbf{U},
\]
where $E (\mathbf{z} \mathbf{z}' \mathbf{z} \mathbf{z}' ) = \sum_{i=1}^p (\kappa
_i+p+2)\mathbf{E}^{ii}$ is diagonal, and therefore the rows of
$\mathbf{U}$ are
the eigenvectors of~$\mathbf{B}$. The order of the eigenvectors is then
given by the order of the corresponding eigenvalues, that is, by the
kurtosis order.
As $\mathbf{W}$ is also affine equivariant, it is an independent component
functional if Assumption~\ref{a3} is replaced
by the following stronger assumption.
%

\begin{assumption}\label{a5}
The fourth moments of $\mathbf{z}$ exist and are distinct.
\end{assumption}
%

\begin{remark}
Notice that Assumption~\ref{a5} $\Rightarrow$ Assumption~\ref{a4}
$\Rightarrow$ Assumption~\ref{a3}.
If Assumption~\ref{a5} is not true and there are only $m<p$ distinct
kurtosis values with
multiplicities $p_1,\ldots,p_m$, FOBI still finds these $m$ subspaces,
and the FOBI solutions at
$\mathbf{z}$ are of the block-diagonal form $\operatorname
{diag}(\mathbf{U}_1,\ldots, \mathbf{U}_m)$
with orthogonal
$p_i\times p_i$ matrices $\mathbf{U}_i$, $i=1,\ldots,m$.
\end{remark}

It is again sufficient to consider the limiting distribution of the estimator
$\hat{\mathbf{W}}=\mathbf{W}(\mathbf{Z})$ only.
Then the asymptotical behavior of the FOBI estimator is given as follows.

%
\begin{theorem}
\label{asymp3}
Let $\mathbf{Z}=(\mathbf{z}_1,\ldots,\mathbf{z}_n)$ be a random
sample from a
distribution of $\mathbf{z}$ with bounded
eighth moments and satisfying the Assumptions \ref{a1}, \ref{a2} and
\ref{a5} with
$\kappa_1 > \cdots> \kappa_p$. Then $\hat{\mathbf{W}}\to_{P}
\mathbf{I}_p$ and
\begin{eqnarray}
&& \sqrt{n} (\hat{w}_{kk}-1) = -\frac{1}2 \sqrt{n} (
\hat{s}_{kk}-1)+o_{P} (1) \quad\mbox{and}
\nonumber
\\
&&\sqrt{n} \hat{w}_{kl}
\nonumber
\\
&&\quad= \biggl(\sqrt{n} \hat{r}_{kl}+ \sqrt{n} \hat{r}_{lk}+
\sqrt{n} \sum_{m\ne k,l}\hat{r}_{mlk}
\nonumber
\\
&&\qquad{}-(\kappa_k+p+4)\sqrt{n} \hat {s}_{kl} \biggr)
\big/({\kappa_k-\kappa_l})+o_{P}(1),
\nonumber
\\
\eqntext{\displaystyle k\neq l.}
\end{eqnarray}
\end{theorem}

For an alternative asymptotic presentation of the $\sqrt{n} \hat w_{kl}$,
see~\citet{Ilmonenetal2010a}. The joint\vspace*{1pt} limiting
multivariate normality of
$\sqrt{n}\cdot\break  \operatorname{vec}(\hat{\mathbf{S}},\operatorname
{off}(\hat{\mathbf{R}}))$ then implies the
following.

%
\begin{corollary}
\label{FOBIasv}
Under the assumptions of Theorem~\ref{asymp3}, the limiting
distribution of
$\sqrt{n} \operatorname{vec} (\hat{\mathbf{W}}- \mathbf{I}_p)$ is a
multivariate normal with zero mean vector and componentwise variances
\begin{eqnarray}
&&\operatorname{ASV}(\hat{w}_{kk}) = (\kappa_k+2)/4\quad
\mbox{and}
\nonumber
\\
&&\operatorname{ASV}(\hat{w}_{kl})
\nonumber
\\
&&\quad= \biggl(\sigma_k^2+\sigma_l^2-
\kappa_k^2-6(\kappa _k+\kappa_l)
\nonumber
\\
&&\hspace*{26pt} {}-22+2p+ \sum_{j\neq k,l}
\kappa_j \biggr) \big/(\kappa_k-\kappa _l)^2,
\nonumber
\\
\eqntext{\displaystyle k\neq l.}
\end{eqnarray}
\end{corollary}

%
\begin{remark}\label{rem7}
Let $\mathbf{x}$ be a $p$-vector with mean vector $\bolds{\mu}$ and
covariance
matrix $\bolds{\Sigma}$.
The FOBI procedure may then be seen also as a comparison of two scatter
functionals, namely,
\begin{eqnarray*}
\operatorname{Cov}(\mathbf{x})&=& \bolds{\Sigma}\quad\mbox{and}
\\
\operatorname{Cov}_4(\mathbf{x}) &=& E \bigl(( \mathbf{x}-\bolds{
\mu}) (\mathbf{x}-\bolds{\mu})' \bolds{\Sigma}^{-1} (
\mathbf{x}-\bolds{\mu}) (\mathbf{x}-\bolds{\mu})' \bigr),
\end{eqnarray*}
and the FOBI functional then satisfies $\mathbf{W} \operatorname
{Cov}(\mathbf{x})\mathbf{W}'=\mathbf{I}_p$ and
$\mathbf{W}\operatorname{Cov}_4(\mathbf{x})\mathbf{W}'\in\mathcal
{D}$. Other
independent component
functionals are obtained if
$\operatorname{Cov}$ and $\operatorname{Cov}_4$ are replaced by any
scatter matrices with the
independence property; see
\citet{Oja2006} and \citet{Tyler2009}.
\end{remark}

\subsection{Joint Use of Kurtosis Matrices $\mathbf{B}^{ij}$---JADE}
\label{sJADE}

The approach in Section~\ref{sFOBI} was based on the fact that the
kurtosis matrix $\mathbf{B}$ is
diagonal at $\mathbf{z}$. As shown before, the fourth cumulant matrices
\[
\mathbf{C}^{ij}=\mathbf{B}^{ij}-\mathbf{E}^{ij}-
\bigl(\mathbf {E}^{ij} \bigr)'-\operatorname{tr} \bigl(
\mathbf{E}^{ij} \bigr)\mathbf{I}_p,\quad i,j=1,\ldots,p,
\]
are also all diagonal at $\mathbf{z}$. Therefore, a natural idea is to try
to find an orthogonal matrix
$\mathbf{U}$ such that the matrices $\mathbf{U}\mathbf
{C}^{ij}\mathbf{U}'$, $i,j=1,\dots
,p$, are all ``as diagonal as
possible.'' In the engineering literature this approach is known as
joint approximate diagonalization
of eigenmatrices (JADE); see \citet{CardosoSouloumiac1993}. The
functional is then defined as
follows.
%

\begin{definition}
The JADE functional is $\mathbf{W}(F_{\mathbf{x}})=\mathbf{U}\bolds
{\Sigma}^{-1/2}$,
where $\bolds{\Sigma}=\operatorname{Cov}(\mathbf{x})$
and the orthogonal matrix $U$ maximizes
\[
\sum_{i=1}^p\sum
_{j=1}^p \bigl\Vert\operatorname{diag} \bigl(\mathbf{U}
\mathbf {C}^{ij} \mathbf{U}' \bigr) \bigr\Vert^2.
\]
\end{definition}

First note that
\begin{eqnarray*}
&& \sum_{i=1}^p\sum
_{j=1}^p \bigl\Vert\operatorname{diag} \bigl(\mathbf{U}
\mathbf {C}^{ij} \mathbf{U}' \bigr) \bigr\Vert^2+
\sum_{i=1}^p \sum
_{j=1}^p \bigl\Vert\operatorname{off} \bigl(\mathbf{U}
\mathbf {C}^{ij} \mathbf{U}' \bigr) \bigr\Vert^2
\\
&&\quad=\sum_{i=1}^p\sum
_{j=1}^p \bigl\Vert\mathbf{C}^{ij}\bigr\Vert^2.
\end{eqnarray*}
The solution thus minimizes the sum of squared off-diagonal elements of
$\mathbf{U} \mathbf{C}^{ij} \mathbf{U}'$,
$i,j=1,\dots,p$. Notice that, at $\mathbf{z}$, the only possible nonzero
elements of $\mathbf{C}^{ij}$,
$i,j=1,\ldots,p$, are $(\mathbf{C}^{ii})_{ii}=\kappa_i$. For the separation
of the components, we therefore need
Assumption~\ref{a4} saying that at most one of the kurtosis values
$\kappa_i$ is zero. The JADE
functional $\mathbf{W}(F)$ is an IC functional, as we can prove in the
following.

%
\begin{theorem}
\label{JADEfunc}
\textup{(i)} Write $\mathbf{x}_{st}=\bolds{\Sigma}^{-1/2}(\mathbf
{x}-\bolds
{\mu})$ for the
standardized random vector from
the IC model~(\ref{ICmodel}) satisfying the Assumptions~\ref{a1},
\ref
{a2}, and \ref{a4}.
If $\mathbf{x}_{st}=\mathbf{U}' \mathbf{z}$, then
\[
D(\mathbf{V})=\sum_{i=1}^p\sum
_{j=1}^p \bigl\Vert\operatorname{diag} \bigl(\mathbf{V}
\mathbf {C}^{ij} \mathbf{V}' \bigr) \bigr\Vert^2,
\quad\mathbf{V}\in \mathcal{U}
\]
is maximized by any $\mathbf{P}\mathbf{J}\mathbf{U}$ where $\mathbf
{P}\in\mathcal{P}$ and
$\mathbf{J}\in\mathcal{J}$.

\textup{(ii)} For any $F_{\mathbf{x}}$ with finite fourth
moments,
$\mathbf{W}(F_{\mathbf{A}\mathbf{x}+\mathbf{b}})=\mathbf{P}\mathbf
{J}\mathbf{W}(F_{\mathbf{x}})\mathbf{A}^{-1}$ for some $\mathbf
{P}\in
\mathcal{P}$ and \mbox{$\mathbf{J}\in\mathcal{J}\!$}.
\end{theorem}

In this case, the matrix $\mathbf{U}=(\mathbf{u}_1,\dots,\mathbf
{u}_p)'$ thus
optimizes the Lagrangian function
\begin{eqnarray*}
L(\mathbf{U},\bolds{\Theta})&=& \sum_{i=1}^p
\sum_{j=1}^p \sum
_{k=1}^p \bigl(\mathbf{u}_k'
\mathbf{C}^{ij} \mathbf{u}_k \bigr)^2-\sum
_{k=1}^p \theta_{kk} \bigl(
\mathbf{u}_k' \mathbf{u}_k-1 \bigr)
\\
&&{}-\sum_{k=1}^{p-1}\sum
_{l=k+1}^{p} \theta_{lk}\mathbf{u}_k'
\mathbf{u}_l,
\end{eqnarray*}
where the symmetric matrix $\bolds{\Theta}=(\theta_{ij})$ contains the
$p(p+1)/2$ Lagrangian multipliers
of the optimization problem. Write
\begin{eqnarray*}
\mathbf{T}(\mathbf{u})&=& \sum_{i=1}^p
\sum_{j=1}^p \bigl(\mathbf{u}'
\mathbf{C}^{ij} \mathbf{u} \bigr)\mathbf {C}^{ij} \mathbf{u}
\quad\mbox{and}
\\
\mathbf{T}( \mathbf{U})&=& \bigl(\mathbf{T}(\mathbf{u}_1),\dots,
\mathbf {T}(\mathbf{u}_p) \bigr)'.
\end{eqnarray*}
The Lagrangian function then yields the estimating equations
\begin{eqnarray*}
\mathbf{u}_i' \mathbf{T}(\mathbf{u}_j) &=&
\mathbf{u}_j' \mathbf {T}(\mathbf{u}_i)
\quad\mbox{and}
\\
\mathbf{u}_i' \mathbf{u}_j &=&
\delta_{ij}, \quad i,j=1,\ldots,p,
\end{eqnarray*}
and the equations suggest a fixed-point algorithm with the steps
$\mathbf{U}\leftarrow\mathbf{T} (\mathbf{T}' \mathbf{T})^{-1/2}$.
The estimating equations can also again be used to find the following
asymptotical distribution
of the JADE estimate $\hat{\mathbf{W}}=\mathbf{W}(\mathbf{Z})$.

%
\begin{theorem}
\label{asymp4}
Let $\mathbf{Z}=(\mathbf{z}_1,\ldots,\mathbf{z}_n)$ be a random
sample from a
distribution of $\mathbf{z}$ with bounded
eighth moments satisfying the Assumptions \ref{a1}, \ref{a2}, and
\ref
{a4}. Then there is a
sequence of solutions $\hat{\mathbf{W}}$ such that $\hat{\mathbf
{W}}\to_{P}
\mathbf{I}_p$ and
\[
\sqrt{n} (\hat{w}_{kk}-1) = -1/2 \sqrt{n} (\hat {s}_{kk}-1)+o_{P}(1),
\quad k=l
\]
and
\begin{eqnarray*}
&&\sqrt{n} \hat{w}_{kl}
\\
&&\quad= \frac{\kappa_{k}\sqrt{n} \hat{r}_{kl}-
\kappa_{l}\sqrt{n} \hat{r}_{lk} + (3\kappa_{l}-3\kappa_{k}-\kappa
_{k}^2 )\sqrt{n} \hat{s}_{kl}}{
\kappa_{k}^2+\kappa_{l}^2}
\\
&&\qquad{} + o_{P}(1), \quad k\neq l.
\end{eqnarray*}
\end{theorem}

%
\begin{corollary}
Under the assumptions of Theorem~\ref{asymp4}, the limiting
distribution of
$\sqrt{n} \operatorname{vec} (\hat{\mathbf{W}}-\mathbf{I}_p)$ is a
multivariate normal
with zero mean vector and
componentwise variances
\begin{eqnarray}
&&\operatorname{ASV}(\hat{w}_{kk}) = (\kappa_k+2)/4 \quad
\mbox{ and}
\nonumber
\\
&&\operatorname{ASV}(\hat{w}_{kl})
\nonumber
\\
&&\quad= \frac{\kappa_k^2(\sigma_k^2-\kappa_k^2-6\kappa
_k-9)+\kappa_l^2(\sigma_l^2
-6\kappa_l-9)}{(\kappa_k^2+\kappa_l^2)^2},
\nonumber
\\
\eqntext{\displaystyle k\neq l.}
\end{eqnarray}
\end{corollary}

%
\begin{remark}
In the literature, there are several alternative algorithms available
for an approximate
diagonalization of several symmetric matrices, but the statistical
properties of the corresponding
estimates are not known. The most popular algorithm is perhaps the
Jacobi rotation algorithm
suggested in~\citet{Clarkson1988}. It appeared in our simulations that
the Jacobi rotation algorithm
is computationally much faster and always provides the same solution as
our fixed-point algorithm.
The limiting distribution with variances and covariances of the
elements of the JADE estimate (but without the standardization step) was
considered also in \citet{BonhommeRobin2009}.
\end{remark}

%
\begin{remark}
The JADE estimate uses $p^2$ fourth moment matrices in order to be
affine equivariant. Therefore, the
computational load of JADE grows quickly with the number of components.
\citet{Miettinenetal2013a}
suggested a quite similar, but faster method, called $k$-JADE. The
$k$-JADE estimate at $F_{\mathbf{x}}$ is
$\mathbf{W}=\mathbf{U}\mathbf{W}_0$, where $\mathbf{W}_0$ is the
FOBI estimate and the
orthogonal matrix $\mathbf{U}$ maximizes
\[
\sum_{|i-j|<k} \bigl\Vert\operatorname{diag} \bigl(\mathbf{U}
\mathbf{C}^{ij} \mathbf{U}' \bigr) \bigr\Vert^2,
\]
where the $\mathbf{C}^{ij}$'s are calculated for $\mathbf
{x}_{st}=\mathbf{W}_0(\mathbf{x}-\bolds{\mu})$.
It seems to us that this estimate is asymptotically equivalent to the
regular JADE
estimate (with much smaller computational load) if the multiplicities
of the distinct
kurtosis values are at most $k$. Detailed studies are, however, still missing.
\end{remark}

\section{Comparison of the Asymptotic Variances of the Estimates}
\label{sare}

First notice that, for all estimates,
$\sqrt{n} (\mathbf{W}(\mathbf{X})-\bolds{\Omega}^{-1})=\sqrt{n}
(\mathbf{W}(\mathbf{Z})-\mathbf{I}_p)\bolds{\Omega}^{-1}$ and the
comparisons can be made using $\hat{\mathbf{W}}=\mathbf{W}(\mathbf
{Z})$ only. Second,
for all estimates,
$
\sqrt{n} (\hat w_{kk}-1)=-1/2 \sqrt{n} (\hat s_{kk}-1)+o_P(1)
$
$k=1,\dots,p$, and therefore the diagonal elements of $\hat{\mathbf{W}}$
should not be used in the
comparison. It is then natural to compare the estimates using the sum
of asymptotic variances of
the off-diagonal elements of $\hat{\mathbf{W}}$, that is,
%
%
\begin{equation}
\label{criterium} \sum_{k=1}^{p-1} \sum
_{l=k+1}^p \bigl( \operatorname{ASV}(
\hat{w}_{kl})+\operatorname {ASV}(\hat {w}_{lk}) \bigr).
\end{equation}
Next note that, for all estimates, except FOBI, the limiting variances of
$\sqrt{n} \hat w_{kl}$, $k\neq l$, surprisingly depend only on the
$k$th and $l$th marginal
distribution (through $\kappa_k$, $\kappa_l$, $\sigma^2_k$, and
$\sigma
^2_l$) and do not depend either
on the number or on the distributions of the other components. Based
on the results in the
earlier sections, we have the following conclusions:
\begin{longlist}[4.]
\item[1.] $\sqrt{n} \hat w_{kl}$ of the symmetric FastICA estimate and
that of the JADE estimate are
asymptotically equivalent, that is, their difference converges to zero
in probability if the $k$th
and $l$th marginal distributions are the same.
\item[2.] If the independent components are identically distributed, then
the symmetric FastICA and
JADE estimates are asymptotically equivalent. In this case, their
criterium value (\ref{criterium})
is one half of that of the deflation-based FastICA estimate. The FOBI
estimate fails in this case.
\item[3.] $\operatorname{ASV}(\hat{w}_{kl})$ of the FOBI estimate is
always larger than
or equal to that for symmetric
FastICA, $k\ne l$. This follows as $\kappa_k\ge-2$ for all $k$. The
larger the other kurtosis
values, the larger is the $\operatorname{ASV}(\hat{w}_{kl})$ of FOBI.
The variances
are equal when $p=2$ and
$\kappa_k>0>\kappa_l$.
\item[4.] $\sqrt{n} \hat w_{kp}$ of the deflation-based FastICA estimate
and of the JADE estimate
are asymptotically equivalent if the $p$th marginal distribution is normal.
\end{longlist}

The criterium value (\ref{criterium}) is thus the sum of the pairwise terms
$\operatorname{ASV}(\hat{w}_{kl})+\operatorname{ASV}(\hat
{w}_{lk})$, which do not depend on the
number or distributions of other
components except for the FOBI estimate. So in most cases the comparison
of the estimates can be made only through the values $\operatorname
{ASV}(\hat
{w}_{kl})+\operatorname{ASV}(\hat{w}_{lk})$. To make
FOBI (roughly) comparable, we use the lower bound of the value
$\operatorname{ASV}(\hat
{w}_{kl})+\operatorname{ASV}(\hat{w}_{lk})$
with $\kappa_j=-2$, $j\ne k,l$; the lower bound is in fact the exact
value in the bivariate case.
In Table~\ref{tab1}, the values $\operatorname{ASV}(\hat
{w}_{kl})+\operatorname{ASV}(\hat{w}_{lk})$
are listed for pairs of
independent components from the following five distributions:
exponential distribution (EX),
logistic distribution (L), uniform distribution (U), exponential power
distribution with shape
parameter value 4 (EP), and normal or Gaussian (G) distribution.
The excess kurtosis values are $\kappa_{\mathrm{EX}}= 6$, $\kappa
_{\mathrm{L}}= 1.2$,
$\kappa_{\mathrm{U}}=- 1.8$,
$\kappa_{\mathrm{EP}}\approx- 0.81$ and $\kappa_{\mathrm{G}}= 0$,
respectively. The
results in Table~\ref{tab1}
are then nicely in accordance with our general notions above and show
that none of the estimates
outperforms all the other estimates.

%
\begin{table}
\caption{The values of $\operatorname{ASV}(\hat
{w}_{kl})+\operatorname{ASV}(\hat{w}_{lk})$ for some
selected $k$th and $l$th
component distributions and for deflation-based FastICA (DFICA),
symmetric FastICA (SFICA),
FOBI, and JADE estimates. For FOBI, the lower bound of
$\operatorname{ASV}(\hat
{w}_{kl})+\operatorname{ASV}(\hat{w}_{lk})$ is
used}\label{tab1}
\begin{tabular*}{\tablewidth}{@{\extracolsep{\fill}}lcccc@{}}
\hline
& \textbf{DFICA} & \textbf{SFICA} & \textbf{FOBI} & \textbf{JADE} \\
\hline
EX--EX & 11.00 & \phantom{0}5.50 & $\infty$ & \phantom{0}5.50 \\
EX--L & 11.00 & \phantom{0}8.52 & 19.18 & 10.22 \\
EX--U & 11.00 & \phantom{0}7.69 & \phantom{0}7.69 & 10.17 \\
EX--EP & 11.00 & \phantom{0}8.63 & \phantom{0}8.63 & 10.61 \\
EX--G & 11.00 & 11.33 & 11.33 & 11.00 \\
L--L & 31.86 & 15.93 & $\infty$ & 15.93 \\
L--U & 31.86 & \phantom{0}8.43 & \phantom{0}8.43 & \phantom{0}8.43 \\
L--EP & 31.86 & 12.38 & 12.38 & 15.63 \\
L--G & 31.86 & 40.19 & 40.19 & 31.86 \\
U--U &  \phantom{0}1.86 & \phantom{0}0.93 & $\infty$ & \phantom{0}0.93 \\
U--EP & \phantom{0}1.86 & \phantom{0}1.80 & 40.63 & \phantom{0}1.50 \\
U--G & \phantom{0}1.86 & 10.19 & 10.19 & \phantom{0}1.86 \\
EP--EP & \phantom{0}6.39 & \phantom{0}3.20 & $\infty$ & \phantom
{0}3.20 \\
EP--G & \phantom{0}6.39 & 34.61 & 34.61 & \phantom{0}6.39 \\
\hline
\end{tabular*}
\end{table}

\begin{figure*}

\includegraphics{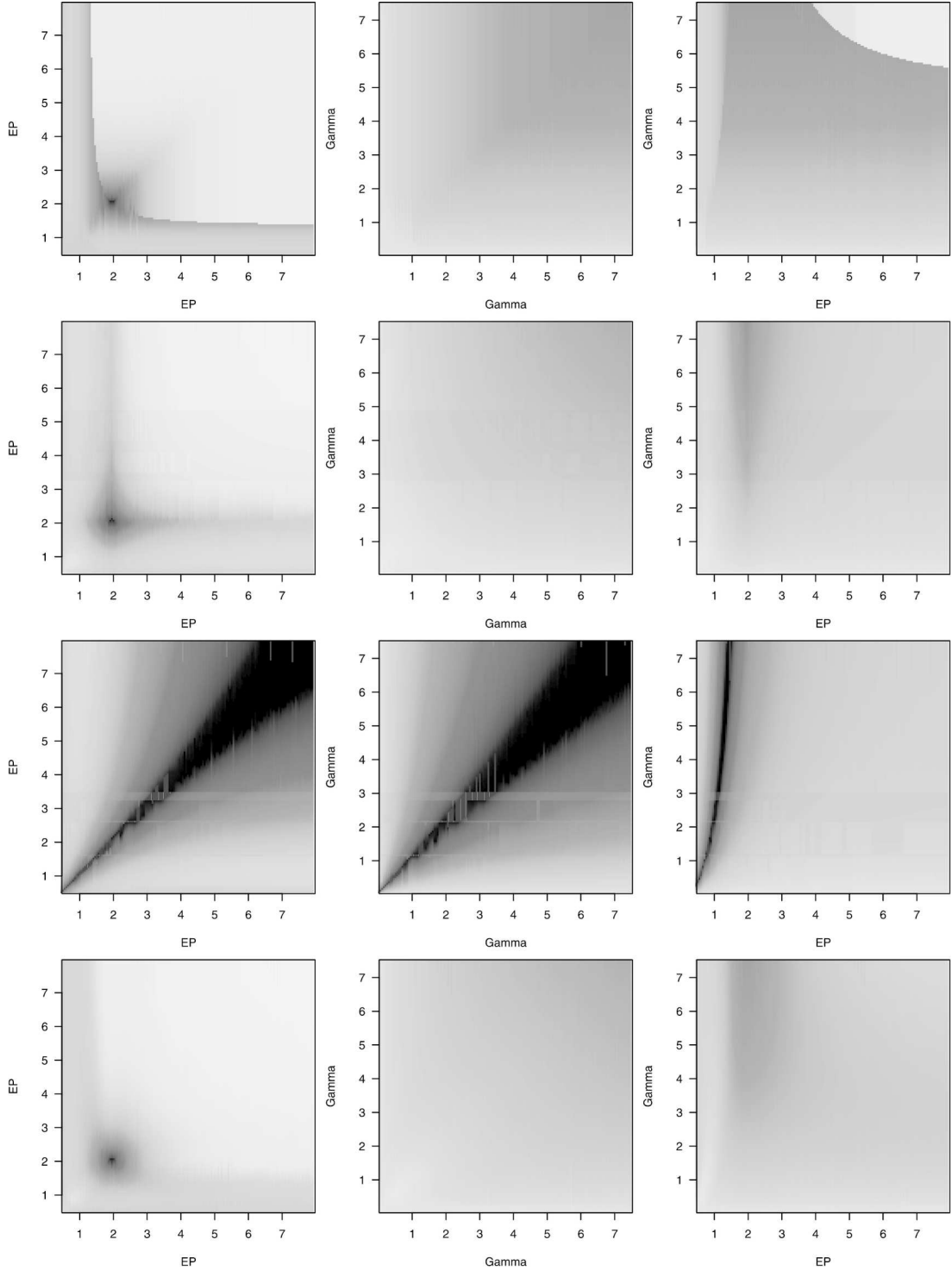}

\caption{Contour maps of $\operatorname{ASV}(\hat
{w}_{kl})+\operatorname{ASV}(\hat{w}_{lk})$ for
different estimates and for
different independent component distributions. The distributions are
either exponential power
distributed (EP) or gamma distributed (Gamma) with varying shape
parameter values. The estimates,
from up to down, are deflation-based FastICA, symmetric FastICA, FOBI,
and JADE. For FOBI, the
lower bound of $\operatorname{ASV}(\hat{w}_{kl})+\operatorname
{ASV}(\hat{w}_{lk})$ is used. The
lighter the color is, the
lower is the variance.}\vspace*{-12pt}
\label{fig1}
\end{figure*}

Further, in Figure~\ref{fig1}, we plot the values $\operatorname
{ASV}(\hat
{w}_{kl})+\operatorname{ASV}(\hat{w}_{lk})$ when the
independent components come (i) from the standardized (symmetric)
exponential power distribution or
(ii) from the standardized (skew) gamma distribution. The limiting
variances then depend only
on the shape parameters of the models. In the plot, the darker the
point, the higher the value and
the worse the estimate. The density function for the exponential power
distribution with zero mean
and variance one and with shape parameter $\beta$ is
\[
f(x)=\frac{\beta\exp\{-(|x|/\alpha)^\beta\}}{2 \alpha\Gamma
(1/\beta)},
\]
where $\beta>0$, $\alpha=(\Gamma(1/\beta)/\Gamma(3/\beta
))^{1/2}$, and
$\Gamma$ is the gamma
function. Notice that $\beta=2$ gives the normal (Gaussian)
distribution, $\beta=1$ gives the
heavy-tailed Laplace distribution, and the density converges to an
extremely low-tailed uniform
density as $\beta\to\infty$. The family of skew distributions for the
variables is coming from
the gamma distribution with shape parameter $\alpha$ and shifted and
rescaled to have mean zero and
variance one. For $\alpha=k/2$, the distribution is a chi-square
distribution with $k$ degrees of
freedom, $k=1,2,\ldots.$ For $\alpha=1$, an exponential distribution is
obtained, and the
distribution is converging to a normal distribution as $\alpha\to
\infty$.

For all estimates, Figure~\ref{fig1} shows that $\operatorname
{ASV}(\hat
{w}_{kl})+\operatorname{ASV}(\hat{w}_{lk})$ gets high
values with $\beta$ close to 2 (normal distribution). Also, the
variances are growing with
increasing $\alpha$. The FOBI estimate is poor if the marginal kurtosis
values are close to each
other. The contours for the deflation-based FastICA estimate illustrate
the fact that the criterium
function $\operatorname{ASV}(\hat{w}_{12})+\operatorname{ASV}(\hat
{w}_{21})$ is not continuous at the
points for which
$\kappa_k+\kappa_l=0$. This is due to the fact that the order in which
the components are found
changes at that point. The symmetric FastICA and JADE estimates are
clearly the best estimates
with minor differences.

%

\section{Discussion}

Many popular methods to solve the independent component analysis problem
are based on the use of
univariate and multivariate fourth moments. Examples include
FOBI~(\citeauthor{Cardoso1989}, \citeyear{Cardoso1989}),
JADE~(\citeauthor{CardosoSouloumiac1993}, \citeyear
{CardosoSouloumiac1993}), and FastICA~(\citeauthor{Hyvarinen1999},
\citeyear{Hyvarinen1999}). In the engineering literature,
these ICA methods have originally been formulated and regarded as
algorithms only, and therefore
the rigorous analysis and comparison of their statistical properties
have been missing until very
recently. The statistical properties of the deflation-based FastICA
method were derived
in~\citet{Ollila2010} and~\citet{Nordhausenetal2011}. The
asymptotical
behavior of the FOBI
estimate was considered in~\citet{Ilmonenetal2010a}, and the
asymptotical distribution of the JADE estimate
(without the standardization step) was considered in~\citet
{BonhommeRobin2009}.
This paper describes in detail the independent
component functionals based on fourth moments through corresponding
optimization problems,
estimating equations, fixed-point algorithms and the assumptions they
need, and provides for the very
first time the limiting statistical properties of the JADE estimate.
Careful comparisons of the
asymptotic variances revealed that, as was expected, JADE and the
symmetric version of FastICA performed
best in most cases. It was surprising, however, that the JADE and
symmetric FastICA estimates are
asymptotically equivalent if the components are identically
distributed. The only noteworthy difference
between these two estimators appeared when one of the components has a
normal distribution. Then JADE
outperforms symmetric FastICA. Recall that JADE requires the
computation of $p^2$ matrices of size
$p\times p$ and, thus, the use of JADE becomes impractical with a large
number of independent components.
On\vadjust{\goodbreak} the other hand, FastICA estimates are
sometimes difficult to find
due to convergence problems of the
algorithms, when the sample size is small.

In this paper we considered only the most basic IC model, where the
number of independent components equals
the observed dimension and where no additive noise is present. In
further research we will consider also
these cases. Note that some properties of JADE for noisy ICA were
considered in \citet{BonhommeRobin2009}.

%
\begin{appendix}
\section*{Appendix: Proofs of the Theorems}

\begin{pf*}{Proof of Theorem~\protect\ref{thm1}}
Let $\bolds{\Omega}=\mathbf{O}\mathbf{D}\mathbf{V}'$ be the
singular value decomposition
of full-rank $\bolds{\Omega}$.
Then $\bolds{\Sigma}=\bolds{\Omega}\bolds{\Omega}'=\mathbf
{O}\mathbf{D}^2\mathbf{O}'$, and
$\bolds{\Sigma}^{-1/2}=\mathbf{O}\mathbf{J}\mathbf{D}^{-1} \mathbf
{O}'$ for some $\mathbf{J}\in
\mathcal{J}$. ($\mathbf{J}$ is needed to
make $\bolds{\Sigma}^{-1/2}$ positive definite.)
Then
\begin{eqnarray*}
\mathbf{x}_{st}&=&\bolds{\Sigma}^{-1/2} (\mathbf{x}-\bolds{
\mu} )=\mathbf{O} \mathbf{J}\mathbf{D}^{-1} \mathbf{O}'
\mathbf{O}\mathbf{D}\mathbf {V}' \mathbf{z}
\\
&=& \mathbf{O}\mathbf{J}\mathbf{V}'\mathbf{z} =\mathbf{U}\mathbf{z}
\end{eqnarray*}
with an orthogonal $\mathbf{U}=\mathbf{O}\mathbf{J}\mathbf{V}'$.
\end{pf*}

\begin{pf*}{Proof of Theorem~\protect\ref{key2}}
If $\mathbf{u}' \mathbf{u}=1$, then it is straightforward to see that
\[
E \bigl[ \bigl(\mathbf{u}'\mathbf{z} \bigr)^4-3 \bigr]=
\sum_{i=1}^p u_i^4
\bigl[E \bigl(z_i^4 \bigr)-3 \bigr].
\]
It then easily follows that
\begin{eqnarray*}
\bigl\vert E \bigl[ \bigl(\mathbf{u}'\mathbf{z} \bigr)^4
\bigr]-3 \bigr\vert &\le& \sum_{i=1}^p
u_i^4 \bigl\vert E \bigl(z_i^4
\bigr)-3 \bigr\vert
\\
&\le&\max_{i=1,\dots,p} \bigl\vert E \bigl(z_i^4
\bigr)-3 \bigr\vert
\end{eqnarray*}
and that, for any orthogonal $\mathbf{U}=(\mathbf{u}_1,\ldots
,\mathbf{u}_p)'$,
\begin{eqnarray*}
\sum_{j=1}^p \bigl\vert E \bigl[ \bigl(
\mathbf{u}_j'\mathbf {z} \bigr)^4 \bigr]-3
\bigr\vert&\le& \sum_{i=1}^p \biggl( \sum
_j u_{ji}^4 \biggr) \bigl\vert E
\bigl(z_i^4 \bigr)-3 \bigr\vert
\\
&\le& \sum_{i=1}^p \bigl\vert E
\bigl(z_i^4 \bigr)-3 \bigr\vert.
\end{eqnarray*}
For the first result, see also Lemma~2 in \citet{BugrienKent2005}.
\end{pf*}

\begin{pf*}{Proof of Theorem~\protect\ref{asymp7}}
As the functions
\begin{eqnarray*}
D_n(\mathbf{U})&=& \sum_{j=1}^p
\Biggl\vert n^{-1} \sum_{i=1}^n \bigl(
\mathbf{u}_j' \mathbf{x}_{st,i}
\bigr)^4 -3 \Biggr\vert\quad\mbox{and}
\\
D(\mathbf{U}) &=& \sum_{j=1}^p \bigl\vert E
\bigl[ \bigl(\mathbf{u}_j' \mathbf{z}
\bigr)^4 \bigr]-3 \bigr\vert
\end{eqnarray*}
are continuous and $D_n(\mathbf{U})\to_P D(\mathbf{U})$ for all
$\mathbf{U}$, then,
due to the compactness of
$\mathcal{U}$,
also
\[
\sup_{\mathbf{U}\in\mathcal{U}} \bigl\vert D_n(\mathbf {U})-D(\mathbf{U})
\bigr\vert\to_P 0.
\]
$D(\mathbf{U})$ attains its maximum at any $\mathbf{J}\mathbf{P}$,
where $\mathbf{J}\in
\mathcal{J}$ and $\mathbf{P}\in\mathcal{P}$.
This further implies that there is a sequence of maximizers that
satisfy $\hat{\mathbf{U}}\to_P \mathbf{I}_p$,
and therefore also $\hat{\mathbf{W}}=\hat{\mathbf{U}} \hat{\mathbf
{S}}^{-1/2}\to_P
\mathbf{I}_p$.

For the estimate $\hat{\mathbf{W}}$, the estimating equations are
\begin{eqnarray*}
\hat{\mathbf{w}}_k' \hat{\mathbf{T}}(\hat{
\mathbf{w}}_l) \hat\pi_l &=& \hat{\mathbf{w}}_l'
\hat{\mathbf{T}}(\hat{ \mathbf{w}}_k) \hat\pi_k \quad
\mbox{and}
\\
\hat{\mathbf{w}}_k' \hat{\mathbf{S}} \hat{
\mathbf{w}}_l &=& \delta_{ij},\quad k,l=1,\ldots,p,
\end{eqnarray*}
where $ \hat{\mathbf{T}}(\hat{\mathbf{w}}_k)=n^{-1} \sum_i (\hat
{\mathbf{w}}_k'
(\mathbf{z}_i-\bar{\mathbf{z}}))^3(\mathbf{z}_i-
\bar{\mathbf{z}})$.
It is straightforward to see that the second set of estimating
equations gives
\renewcommand{\theequation}{\arabic{equation}}
\setcounter{equation}{4}
%
\begin{eqnarray}
\sqrt{n} (\hat{w}_{kk}-1)&=& -2^{-1} \sqrt{n} (
\hat{s}_{kk}-1)+o_{P}(1)\quad\mbox{and}\hspace*{-12pt}
\nonumber
\\[-8pt]
\label{wlk}
\\[-8pt]
\nonumber
\sqrt{n} (\hat{w}_{kl}+\hat{w}_{lk})&=& -\sqrt{n}
\hat{s}_{kl}+o_{P}(1).
\end{eqnarray}

Consider then the first set of estimating equations for $k\neq l$. To
shorten the notation, write
$\hat{\mathbf{T}}(\hat{\mathbf{w}}_k)=\hat{\mathbf{T}}_k$. Now
\[
\sqrt{n} \hat{\mathbf{w}}_k'\hat{\mathbf{T}}_l=
\sqrt{n} (\hat{\mathbf{w}}_k-\mathbf{e}_k)'
\hat{\mathbf{T}}_l+ \sqrt{n} \mathbf{e}_k'(
\hat{\mathbf{T}}_l-\beta_l \mathbf{e}_l).
\]
Using equation (2) in \citet{Nordhausenetal2011} and Slutsky's
theorem, the above equation reduces to
\begin{eqnarray*}
&&\sqrt{n} \hat{\mathbf{w}}_k'\hat{
\mathbf{T}}_l \hat \pi_l
\\
&&\quad= \bigl(\sqrt{n} (\hat{\mathbf{w}}_k-\mathbf{e}_k)'
\beta _l\mathbf{e}_l + \mathbf{e}_k'
\bigl(\sqrt{n} {\hat{\mathbf{T}}}^*_l- \gamma_l
\mathbf{e}_l \mathbf{e}_l' \sqrt{n} \bar{
\mathbf{x}}
\\
&&\hspace*{109pt}\qquad{}+\bolds{\Delta}_l\sqrt{n} (\hat{\mathbf
{w}}_l-\mathbf{e}_l) \bigr) \bigr)\pi_l
\\
&&\qquad{}+o_{P}(1),
\end{eqnarray*}
where $\hat{\mathbf{T}}^*_l=n^{-1} \sum_i((\mathbf{e}_l' \mathbf
{z}_i)^3-\gamma
_l)\mathbf{z}_i$ and
$\bolds{\Delta}_l=\break 3 E[(\mathbf{e}_l' \mathbf{z}_i)^2\mathbf
{z}_i\mathbf{z}_i' ]$.
According to our estimating equation, the above expression should be
equivalent to
\begin{eqnarray*}
\sqrt{n} \hat{\mathbf{w}}_l' \hat{
\mathbf{T}}_k \hat\pi_k &=& \bigl(\sqrt{n} (\hat{
\mathbf{w}}_l-\mathbf {e}_l)'
\beta_k\mathbf{e}_k
\\
&&\hspace*{3pt} {}+ \mathbf{e}_l' \bigl(\sqrt{n} \hat{
\mathbf{T}}^*_k- \gamma_k\mathbf{e}_k
\mathbf{e}_k'\sqrt{n} \bar{\mathbf{z}}
\\
&&\hspace*{3pt} {}+\bolds{\Delta}_k\sqrt{n} (\hat{
\mathbf{w}}_i- \mathbf{e}_i) \bigr) \bigr)
\pi_k+o_{P}(1).
\end{eqnarray*}
This further implies that
\begin{eqnarray*}
&& (\beta_l\pi_l-3\pi_k)\sqrt{n}
\hat{w}_{kl}-(\beta_k\pi_k-3\pi_l)
\sqrt{n} \hat{w}_{lk}
\\
&&\quad=\sqrt{n} (\hat r_{kl} \pi_k+\hat r_{lk}
\pi_l)+o_{P}(1),
\end{eqnarray*}
where $\hat r_{kl}=\sum_i(z_{ik}^3-\gamma_k)z_{il}$. Now using~(\ref
{wlk}), we have that
\begin{eqnarray*}
&&(\beta_l\pi_l-3\pi_k)\sqrt{n}
\hat{w}_{kl}
\\
&&\qquad{}+(\beta_k\pi_k-3\pi_l) (\sqrt{n}
\hat{s}_{kl}+\sqrt{n} \hat{w}_{kl})
\\
&&\quad=\sqrt{n} (\hat r_{kl}\pi_k+ \hat r_{lk}
\pi_l)+o_{P}(1).
\end{eqnarray*}
Then
\begin{eqnarray*}
&& \bigl(|\beta_k-3|+|\beta_l-3| \bigr)\sqrt{n}
\hat{w}_{kl}
\\
&&\quad=\sqrt{n} (\hat r_{kl}\pi_k+\hat r_{lk}
\pi_l)+(3\pi_l-\beta_k\pi _k)
\sqrt{n} \hat{s}_{kl}
\\
&&\qquad{}+o_{P}(1),
\end{eqnarray*}
which gives the desired result.
\end{pf*}

\begin{pf*}{Proof of Theorem~\ref{asymp3}}
As mentioned in Remark~\ref{rem7}, the FOBI functional diagonalizes
the scatter matrices
$\operatorname{Cov}(\mathbf{x})$ and
$\operatorname{Cov}_4(\mathbf{x})=E[(\mathbf{x}-\bolds{\mu
})(\mathbf{x}-\bolds{\mu
})' \bolds{\Sigma}^{-1}(\mathbf{x}-\bolds{\mu})(\mathbf{x}-\bolds
{\mu})']$ simultaneously.
Then\break $\operatorname{Cov}(\mathbf{z})=\mathbf{I}_p$ and
$\operatorname
{Cov}_4(\mathbf{z})=\mathbf{D}$ with strictly decreasing
diagonal elements. Next write $\hat{\mathbf{S}}$ and $\hat{\mathbf
{S}}_4$ for the
empirical scatter matrices.
Then $\hat{\mathbf{S}}\to_P \mathbf{I}_p$ and $\hat{\mathbf
{S}}_4\to_P \mathbf{D}$,
and, as
$\hat{\mathbf{W}}$ is a continuous
function of $(\hat{\mathbf{S}},\hat{\mathbf{S}}_4)$ in a
neighborhood of $(\mathbf{I}_p,\mathbf{D})$, also
$\hat{\mathbf{W}}\to_P \mathbf{I}_p$.

Let $\tilde{\mathbf{Z}}=(\tilde{\mathbf{z}}_1,\dots,\tilde
{\mathbf{z}}_n)=(\mathbf{z}_1-\bar
{\mathbf{z}},\dots,\mathbf{z}_n-
\bar{\mathbf{z}})$ denote the centered sample,
\begin{eqnarray*}
\sqrt{n} (\hat{\mathbf{S}}_4-\mathbf{D})&=& n^{-1/2}\sum
_{i=1}^n \bigl(\tilde{\mathbf{z}}_i
\tilde{ \mathbf{z}}_i' \hat{\mathbf{S}}^{-1}
\tilde{\mathbf{z}}_i \tilde {\mathbf{z}}_i'-
\mathbf{D} \bigr)
\\
&=& -n^{-1}\sum_{i=1}^n \tilde{
\mathbf{z}}_i\tilde{\mathbf{z}}_i'\sqrt{n}
(\hat{ \mathbf{S}}-\mathbf{I}_p) \tilde{\mathbf{z}}_i
\tilde{\mathbf {z}}_i'
\\
&&{}+n^{-1/2}\sum_{i=1}^n \bigl(
\tilde{\mathbf{z}}_i\tilde{ \mathbf{z}}_i'
\tilde{\mathbf{z}}_i\tilde{\mathbf{z}}_i'-
\mathbf {D} \bigr),
\end{eqnarray*}
where the $(k,l)$ element, $k\neq l$, of the first matrix is
\begin{eqnarray*}
&&-n^{-1}\sum_{i=1}^n
\tilde{z}_{ki}\tilde{z}_i'\sqrt{n} (\hat{
\mathbf{S}}-\mathbf {I}_p)\tilde{z}_i
\tilde{z}_{li}
\\
&&\quad= -2n^{-1}\sum_{i=1}^n
\tilde{z}_{ki}^2\tilde{z}_{li}^2
\sqrt{n} \hat{s}_{kl}+o_P(1)
\\
&&\quad= -2\sqrt{n} \hat{s}_{kl}+o_P(1),
\end{eqnarray*}
and the $(k,l)$ element of the second matrix is
\begin{eqnarray*}
&& n^{-1/2}\sum_{i=1}^n
\tilde{z}_{ki}^3\tilde{z}_{li}+n^{-1/2}
\sum_{i=1}^n \tilde{z}_{ki}
\tilde{z}_{li}^3
\\
&&\quad{}+n^{-1/2}\sum_{i=1}^n
\sum_{m\neq
k,l} \tilde{z}_{mi}^2
\tilde{z}_{ki} \tilde{z}_{li}.
\end{eqnarray*}
Thus,
\[
\sqrt{n} (\hat{\mathbf{S}}_4)_{kl}=\sqrt{n}
\hat{r}_{kl}+\sqrt{n} \hat{r}_{lk}+\sum
_{m\neq k,l} \hat{r}_{mkl}+o_P(1).
\]

Then Theorem~3.1 of \citet{Ilmonenetal2010a} gives
\begin{eqnarray*}
&&\sqrt{n} \hat{w}_{kl}
\\[-2pt]
&&\quad=\frac{\sqrt{n} (\hat{\mathbf{S}}_4)_{kl}-(\kappa
_k+p+2)\sqrt{n} \hat{s}_{kl}}{
\kappa_k+p+2-(\kappa_l+p+2)}+o_P(1)
\\[-2pt]
&&\quad= \biggl(\sqrt{n} \hat{r}_{kl}+\sqrt{n} \hat{r}_{lk}+
\sum_{m\neq k,l} \sqrt{n} \hat{r}_{mkl}
\\[-2pt]
&&\qquad\hspace*{6pt} {}- (\kappa_k+p+4)\sqrt{n}
\hat{s}_{kl} \biggr) \Big/(\kappa_k-\kappa _l)+o_P(1).
\end{eqnarray*}
\vspace*{2pt}\upqed
\end{pf*}

To prove Theorem~\ref{JADEfunc}, we need the following lemma.

%
\begin{lemma}
\label{add}
Denote
\begin{eqnarray*}
\mathbf{C}(\mathbf{x},\mathbf{A}) &=& E \bigl[ \bigl(\mathbf{x}'
\mathbf{A} \mathbf{x} \bigr)\mathbf{x} \mathbf{x}' \bigr]-\mathbf{A}
- \mathbf{A}' -\operatorname {tr}(\mathbf{A})\mathbf{I}_p,
\\[-1pt]
\mathbf{C}^{ij}(\mathbf{x}) &=& E \bigl[ \bigl(\mathbf{x}'
\mathbf {E}^{ij} \mathbf{x} \bigr)\mathbf{x} \mathbf{x}'
\bigr]- \mathbf {E}^{ij} - \mathbf{E}^{ji} -\operatorname{tr}
\bigl(\mathbf {E}^{ij} \bigr)\mathbf{I}_p,
\end{eqnarray*}
where $\mathbf{E}^{ij}=\mathbf{e}_i\mathbf{e}_j'$, $i,j=1,\ldots
,p$. Then $\mathbf{C}(\mathbf{x},\mathbf{A})$ is additive in
$\mathbf{A}=(a_{ij})$,\vspace*{-2pt} that is,
\[
\mathbf{C}(\mathbf{x},\mathbf{A})= \sum_{i=1}^p
\sum_{j=1}^p \mathbf{a}_{ij}
\mathbf{C}^{ij}(\mathbf{x}).
\]
Also, for an orthogonal $\mathbf{U}$, it holds that
\[
\mathbf{C}(\mathbf{U}\mathbf{x},\mathbf{A})=\mathbf{U} \mathbf {C} \bigl(
\mathbf{x},\mathbf{U}' \mathbf{A} \mathbf{U} \bigr)
\mathbf{U}'.
\]
\end{lemma}

\begin{pf} 
For additivity, it is straightforward to see that, for all $\mathbf
{A}, \mathbf{A}_1, \mathbf{A}_2$, and $b$,
\begin{eqnarray*}
\mathbf{C}(\mathbf{x},b \mathbf{A})& =& b \mathbf{C}(\mathbf {x},\mathbf{A})
\quad\mbox{and}
\\
\mathbf{C}(\mathbf{x},\mathbf {A}_1+\mathbf{A}_2) &=&
\mathbf{C}(\mathbf{x},\mathbf{A}_1)+\mathbf{C}( \mathbf{x},
\mathbf{A}_2).
\end{eqnarray*}
For orthogonal $\mathbf{U}$, we obtain
\begin{eqnarray*}
&& \mathbf{C}(\mathbf{U}\mathbf{x},\mathbf{A})
\\
&&\quad= E \bigl[ \bigl(\mathbf{x}'\mathbf{U}'
\mathbf{A}\mathbf {U}\mathbf{x} \bigr) \bigl(\mathbf{U}\mathbf{x}
\mathbf{x}' \mathbf{U}' \bigr) \bigr] -\mathbf{A}-
\mathbf{A}'-\operatorname{tr}(\mathbf{A})\mathbf{I}_p
\\
&&\quad= \mathbf{U} \bigl( E \bigl[ \bigl(\mathbf{x}' \bigl(
\mathbf{U}' \mathbf{A}\mathbf{U} \bigr)\mathbf{x} \bigr)\mathbf{x}
\mathbf{x}' \bigr]- \bigl(\mathbf{U}' \mathbf{A}
\mathbf{U} \bigr)- \bigl(\mathbf{U}' \mathbf{A}\mathbf{U} \bigr)
\\
&&\qquad{}-\operatorname{tr} \bigl( \bigl(\mathbf{U}'\mathbf{A}
\mathbf{U} \bigr) \bigr)\mathbf {I}_p \bigr) \mathbf{U}'
\\
&&\quad= \mathbf{U} \mathbf{C} \bigl(\mathbf{x},\mathbf{U}'
\mathbf{A}\mathbf{U} \bigr) \mathbf{U}'.
\end{eqnarray*}
\upqed
\end{pf}

\begin{pf*}{Proof of Theorem~\protect\ref{JADEfunc}}
(i) First notice that
\begin{eqnarray*}
\mathbf{C}^{ij}(\mathbf{z}) &=& \mathbf{0},\quad\mbox{for $i,j=1,
\ldots,p$ and $i\neq j$}
\\
\mathbf{C}^{ii}(\mathbf{z}) &=& \kappa_i
\mathbf{E}^{ii},\quad \mbox{for $i=1,\ldots,p$}.
\end{eqnarray*}
It then follows that, for an orthogonal $\mathbf{U}=(\mathbf
{u}_1,\ldots,\break 
\mathbf{u}_p)$,\vspace*{-1pt}
\begin{eqnarray*}
C^{ij} \bigl(\mathbf{U}' \mathbf{z} \bigr) &=&
\mathbf{U}' \mathbf {C} \bigl(\mathbf{z},\mathbf{U}
\mathbf{E}^{ij} \mathbf{U}' \bigr)\mathbf{U}
\\
& =& \mathbf{U}' \mathbf{C} \Biggl(\mathbf{z},\sum
_{k=1}^p \sum_{l=1}^p
u_{ki}u_{lj} \mathbf{E}_{kl} \Biggr) \mathbf{U}
\\
&=& \mathbf{U}' \Biggl(\sum_{k=1}^p
\sum_{l=1}^p u_{ki}u_{lj}
\mathbf{C}(\mathbf{z},\mathbf{E}_{kl}) \Biggr) \mathbf{U}
\\
&=& \mathbf{U}' \Biggl(\sum_{k=1}^p
\kappa_k u_{ki}u_{kj} \mathbf{E}_{kk}
\Biggr) \mathbf{U}.
\end{eqnarray*}
Now
\begin{eqnarray*}
D(\mathbf{V})&=&\sum_{i=1}^p \sum
_{j=1}^p \bigl\Vert\operatorname{diag} \bigl(
\mathbf{V}\mathbf{C}^{ij} \bigl(\mathbf{U}' \mathbf{z}
\bigr)\mathbf{V}' \bigr) \bigr\Vert^2
\\
&=& \sum_{i=1}^p \sum
_{j=1}^p \Biggl\Vert\mathbf{V}\mathbf{U}'
\Biggl(\sum_{k=1}^p \kappa_k
u_{ki}u_{kj}E_{kk} \Biggr) \bigl(\mathbf{V}
\mathbf{U}' \bigr)' \Biggr\Vert^2.
\end{eqnarray*}
If we write $\mathbf{G} =\mathbf{V}\mathbf{U}' =(\mathbf
{g}_1,\ldots,\mathbf{g}_p)$, then
$D(\mathbf{V})$ simplifies to
\begin{eqnarray*}
D(\mathbf{V}) &=& \sum_{i=1}^p \sum
_{j=1}^p \sum_{k=1}^p
\Biggl(\sum_{l=1}^p g_{kl}^2
\kappa_l u_{li}u_{lj} \Biggr)^2
\\
&=& \sum_{i=1}^p \sum
_{j=1}^p \sum_{k=1}^p
\sum_{l=1}^p \sum
_{l^*=1}^p g_{kl}^2g_{kl^*}^2
\kappa_l \kappa_{l^*} u_{li}u_{lj}
u_{l^*i}u_{l^*j}
\\
&=& \sum_{k=1}^p \sum
_{l=1}^p \sum_{l^*=1}^p
g_{kl}^2g_{kl^*}^2
\kappa_l \kappa_{l^*}
\\
&&{}\cdot\sum_{i=1}^p
(u_{li}u_{l^*i})\sum_{j=1}^p
(u_{lj} u_{l^*j})
\\
&=& \sum_{k=1}^p \sum
_{l=1}^p g_{kl}^4
\kappa_l^2,
\end{eqnarray*}
which is maximized by $\mathbf{V}=\mathbf{P}\mathbf{J}\mathbf{U}$
for any $\mathbf{P}\in
\mathcal
{P}$ and
$\mathbf{J}\in\mathcal{J}$.

(ii) Write $\mathbf{y}=\mathbf{A}\mathbf{x}+\mathbf{b}$, and let
$\bolds{\mu}$ and
$\bolds{\Sigma}$
denote the mean vector and
covariance matrix of $\mathbf{x}$, respectively. As
\[
\bigl(\mathbf{A}\bolds{\Sigma}\mathbf{A}' \bigr)^{-1/2}
\bigl(\mathbf{A} \bolds{\Sigma}\mathbf{A}' \bigr) \bigl(\mathbf{A}
\bolds{\Sigma }\mathbf{A}' \bigr)^{-1/2}=
\mathbf{I}_p,
\]
we have that
$(\mathbf{A}\bolds{\Sigma}\mathbf{A}')^{-1/2}\mathbf{A}=\mathbf
{Q} \bolds{\Sigma
}^{-1/2}$ for some
$\mathbf{Q}\in\mathcal{U}$,
and therefore $\mathbf{y}_{st}=\mathbf{Q}\mathbf{x}_{st}$ with the
same $\mathbf{Q}\in
\mathcal{U}$.

We thus define $\mathbf{W}(F_{\mathbf{x}})=\mathbf{U}\bolds{\Sigma
}^{-1/2}$, where
$\mathbf{U}$
maximizes the function
\[
D_{\mathbf{x}_{st}}(\mathbf{V})=\sum_{i=1}^p
\sum_{j=1}^p \bigl\Vert\operatorname{diag}
\bigl(\mathbf{V}\mathbf {C}^{ij}(\mathbf{x}_{st})
\mathbf{V}' \bigr) \bigr\Vert^2.
\]
The maximizer $\mathbf{U}$ is not unique, as the maximum is then attained
for any $\mathbf{P}\mathbf{J}\mathbf{U}$ where
$\mathbf{P}\in\mathcal{P}$ and $\mathbf{J}\in\mathcal{J}$.

Consider next the criterium function for the standardized transformed
random variable $\mathbf{y}_{st}$. Then
\begin{eqnarray*}
D_{y_{st}}(\mathbf{V})&=&\sum_{i=1}^p
\sum_{j=1}^p \bigl\Vert\operatorname{diag}
\bigl(\mathbf{V} \mathbf{C}^{ij}(\mathbf{Q} \mathbf{x}_{st})
\mathbf{V}' \bigr) \bigr\Vert^2
\\
&=&\sum_{i=1}^p \sum
_{j=1}^p \bigl\Vert\mathbf{V}\mathbf{Q} \mathbf{C} \bigl(
\mathbf {x}_{st},\mathbf{Q}' \mathbf{E}^{ij}
\mathbf{Q} \bigr) \mathbf {Q}'\mathbf{V}'
\bigr\Vert^2.
\end{eqnarray*}
If we write $\mathbf{G} = \mathbf{V}\mathbf{Q} = (\mathbf
{g}_1,\dots,\mathbf{g}_p)'$, then
\begin{eqnarray*}
&& D_{y_{st}}(\mathbf{V})
\\
&&\quad= \sum_{i=1}^p \sum
_{j=1}^p \sum_{k=1}^p
\bigl(\mathbf{g}_k \mathbf{C} \bigl(\mathbf {x}_{st},
\mathbf{Q}' \mathbf{E}^{ij}\mathbf{Q} \bigr)
\mathbf{g}_k' \bigr)^2
\\
&&\quad=\sum_{i=1}^p \sum
_{j=1}^p \sum_{k=1}^p
\Biggl(\sum_{l=1}^p \sum
_{m=1}^p \sum_{s=1}^p
\sum_{t=1}^p g_{ks}g_{kt}q_{il}q_{jm}
\\
&&\hspace*{140pt} {}\cdot\mathbf{C} \bigl(\mathbf{x}_{st},\mathbf
{E}^{lm} \bigr)_{st} \Biggr)^2
\\
&&\quad=\sum_{i,j,k,l,l^*,m,m^*,s,s^*,t,t^*=1}^p
g_{ks}g_{kt}g_{ks^*}g_{kt^*}q_{il}
\\
&& \hspace*{95pt}\qquad{}\cdot q_{jm}q_{il^*}q_{jm^*}
\mathbf{C} \bigl(\mathbf{x}_{st},\mathbf{E}^{lm}
\bigr)_{st}
\\
&&\hspace*{95pt}\qquad{}\cdot\mathbf{C} \bigl(\mathbf {x}_{st},
\mathbf{E}^{l^*m^*} \bigr)_{s^*t^*}
\\
&&\quad=\sum_{k,l,l^*,m,m^*,s,s^*,t,t^*=1}^p
g_{ks}g_{kt}g_{ks^*}g_{kt^*} \mathbf{C}
\bigl(\mathbf{x}_{st},\mathbf{E}^{lm} \bigr)_{st}
\\
&&\hspace*{83pt}\qquad{}\cdot\mathbf{C} \bigl(\mathbf{x}_{st},
\mathbf{E}^{l^*m^*} \bigr)_{s^*t^*}
\\
&&\qquad{}\cdot\sum_{i}(u_{il}u_{il^*})
\sum_{j}(u_{jm}u_{jm^*})
\\
&&\quad=\sum_{k,l,m,s,s^*,t,t^*=1}^p
g_{ks}g_{kt}g_{ks^*}g_{kt^*} \mathbf{C}
\bigl(\mathbf{x}_{st},\mathbf{E}^{lm} \bigr)_{st}
\\
&&\hspace*{60pt}\qquad{}\cdot\mathbf{C} \bigl(\mathbf{x}_{st},
\mathbf{E}^{lm} \bigr)_{s^*t^*}
\\
&&\quad=D_{\mathbf{x}_{st}}( \mathbf{G}).
\end{eqnarray*}
Hence, $D_{y_{st}}(\mathbf{V})=D_{x_{st}}(\mathbf{V}\mathbf{Q})\leq
D_{\mathbf{x}_{st}}(\mathbf{U})$, with equality, if
$\mathbf{V}=\mathbf{P}\mathbf{J}\mathbf{U}\mathbf{Q}' $ for any
$\mathbf{P}\in\mathcal{P}$ and $\mathbf{J}\in\mathcal{J}$.
Thus, $W(F_{\mathbf{y}})=\mathbf{P}\mathbf{J}\mathbf{U}\mathbf{Q}'
\mathbf{Q}\bolds{\Sigma
}^{-1/2}\mathbf{A}^{-1}=\mathbf{P}\mathbf{J}\mathbf{U}\bolds
{\Sigma}^{-1/2}
\mathbf{A}^{-1}$ for any $\mathbf{P}\in\mathcal{P}$ and $\mathbf
{J}\in\mathcal{J}$.
\end{pf*}

For Theorem~\ref{asymp4} we need the following lemma.
%

\begin{lemma}
\label{JointDiag}
Assume that $\hat{\mathbf{S}}_k$, $k=1,\ldots,K$ are $p\times p$ matrices
such that
$\sqrt{n} (\hat{\mathbf{S}}_k-\bolds{\Lambda}_k)$ are
asymptotically normal
with mean zero and
$\bolds{\Lambda}_k=\operatorname{diag}(\lambda_{k1},\ldots,
\lambda_{kp})$.
Let $\hat{\mathbf{U}}=(\hat{\mathbf{u}}_1,\ldots,\hat{\mathbf
{u}}_p)$ be the orthogonal
matrix that maximizes
\[
\sum_{k=1}^K \bigl\| \operatorname{diag}
\bigl(\hat{\mathbf{U}}' \hat{\mathbf{S}}_k \hat{
\mathbf{U}} \bigr) \bigr\|^2.
\]
Then
\[
\sqrt{n} \hat{u}_{ij}=\frac{\sum_{k=1}^K (\lambda_{ki}-\lambda
_{kj})\sqrt{n} (\hat{\mathbf{S}}_k)_{ij}}{
\sum_{k=1}^K (\lambda_{ki}-\lambda_{kj})^2}+o_{P}(1).
\]
\end{lemma}

\begin{pf} 
The proof is similar to the proof of Theorem~4.1 of \citet
{Miettinenetal2014}.
\end{pf}

\begin{pf*}{Proof of Theorem~\ref{asymp4}}
As the criterium functions
\begin{eqnarray*}
D_n(\mathbf{U})&=& \sum_{j=1}^p
\sum_{j=1}^p \bigl\Vert\operatorname{diag}
\bigl(\mathbf{U}\hat {\mathbf{C}}^{ij} \mathbf{U}' \bigr)
\bigr\Vert^2 \quad\mbox{and}
\\
D( \mathbf{U}) &=& \sum_{j=1}^p \sum
_{j=1}^p \bigl\Vert\operatorname{diag} \bigl(
\mathbf{U} \mathbf {C}^{ij} \mathbf{U}' \bigr)
\bigr\Vert^2
\end{eqnarray*}
are continuous and $D_n(\mathbf{U})\to_P D(\mathbf{U})$ for all
$\mathbf{U}$, then,
due to the compactness of
$\mathcal{U}$,
\[
\sup_{\mathbf{U}\in\mathcal{U}} \bigl\vert D_n(\mathbf {U})-D(\mathbf{U})
\bigr\vert\to_P 0.
\]
$D(\mathbf{U})$ attains its maximum at any $\mathbf{J}\mathbf{P}$
where $\mathbf{J}\in
\mathcal{J}$ and
$\mathbf{P}\in\mathcal{P}$. This further implies that there is a
sequence of
maximizers that satisfy
$\hat{\mathbf{U}}\to_P \mathbf{I}_p$, and therefore also
$\hat{\mathbf{W}}=\hat{\mathbf{U}} \hat{\mathbf{S}}^{-1/2}\to_P
\mathbf{I}_p$.

Let $\tilde{\mathbf{Z}}=(\tilde{\mathbf{z}}_1,\dots,\tilde
{\mathbf{z}}_n)=(\mathbf{z}_1-\bar
{\mathbf{z}},\dots,\mathbf{z}_n-\bar{\mathbf{z}})$
denote the centered sample, and write
\begin{eqnarray*}
\label{T} &&\hat{\bolds{\beta}}=\bolds{\beta}(\tilde{\mathbf{Z}})=n^{-1}
\sum_{i=1}^n \bigl(\tilde{
\mathbf{z}}_i \tilde{\mathbf {z}}_i' \bigr)
\otimes \bigl( \tilde{\mathbf{z}}_i\tilde{\mathbf{z}}_i'
\bigr).
\end{eqnarray*}
As the eighth moments of $\mathbf{z}$ exist, $\sqrt{n} (\hat{\bolds
{\beta}}-
\bolds{\beta})$ is asymptotically
normal with the expected value zero, and $\bolds{\beta}$ as in
Theorem~\ref{B}.

Consider first a general sample whitening matrix $\hat{\mathbf{V}}$ satisfying
$\sqrt{n} (\hat{\mathbf{V}}-\mathbf{I}_p)=O_P(1)$. For the whitened
data we obtain
\begin{eqnarray*}
\tilde{\bolds{\beta}}&=&\bolds{\beta}(\hat{\mathbf{V}} \tilde {
\mathbf{Z}})=n^{-1} \sum_{i=1}^n
\bigl(\hat{\mathbf{V}}\tilde{\mathbf{z}}_i\tilde {
\mathbf{z}}_i' \hat{\mathbf{V}}' \bigr)
\otimes \bigl(\hat{\mathbf {V}} \tilde{\mathbf{z}}_i \tilde{
\mathbf{z}}_i' \hat{\mathbf{V}}' \bigr)
\\
&=& (\hat{ \mathbf{V}}\otimes\hat{\mathbf{V}}) \hat{\bolds{\beta}} \bigl(\hat{
\mathbf{V}}' \otimes\hat{\mathbf{V}}' \bigr),
\end{eqnarray*}
and, further,
\begin{eqnarray*}
&& \sqrt{n} (\tilde{\bolds{\beta}}-\bolds{\beta})
\\
&&\quad= \sqrt{n} (\hat{ \bolds{\beta} }-\bolds{\beta})
\\
&&\qquad{}+ \bigl[ \bigl(\sqrt{n} (\hat{\mathbf{V}}- \mathbf {I}_p)
\otimes\mathbf{I}_p \bigr)+ \bigl( \mathbf{I}_p \otimes
\sqrt{n} (\hat{\mathbf{V}}- \mathbf{I}_p) \bigr) \bigr]\bolds{\beta}
\\
&&\qquad{}+ \bolds{\beta} \bigl[ \bigl(\sqrt{n} \bigl(\hat {
\mathbf{V}}'-\mathbf{I}_p \bigr)\otimes
\mathbf{I}_p \bigr)
\\
&&\hspace*{23pt}\qquad{}+ \bigl(\mathbf{I}_p \otimes\sqrt{n} \bigl(
\hat{\mathbf{V}}'- \mathbf{I}_p \bigr) \bigr) \bigr].
\end{eqnarray*}
Write next
\begin{eqnarray*}
\hat{\mathbf{B}}^{kl}&=& \mathbf{B} \bigl(\mathbf{E}^{kl},
\tilde {\mathbf{Z}} \bigr)=n^{-1} \sum_{i=1}^n
\bigl( \tilde{\mathbf{z}}_i \tilde{\mathbf{z}}_i'
\mathbf{E}^{kl} \tilde {\mathbf{z}}_i \tilde{
\mathbf{z}}_i' \bigr) \quad\mbox{and}
\\
\hat{\mathbf{T}}^{kl} &=& \operatorname{vec} \bigl(\hat{\mathbf
{B}}^{kl} \bigr)=\hat{\bolds{\beta}} \operatorname{vec} \bigl(
\mathbf{E}^{kl} \bigr).
\end{eqnarray*}
Then $\sqrt{n} (\hat{\mathbf{T}}^{kl}-\operatorname{vec}(\mathbf
{B}^{kl}))$ is
asymptotically normal with expected
value zero and $\mathbf{B}^{kl}$ as given in Theorem~\ref{B}. Also,
\[
\sqrt{n} \hat{b}_{kl}^{kk}=\sqrt{n} \bigl(\hat{\mathbf
{B}}^{kk} \bigr)_{kl}=\sqrt{n} \hat{r}_{kl}+o_P(1).
\]
Next, let
\[
\tilde{\mathbf{B}}^{kl}=\mathbf{B} \bigl(\mathbf{E}^{kl},
\hat {\mathbf{V}} \tilde{\mathbf{Z}} \bigr)=n^{-1} \sum
_{i=1}^n \bigl(\hat{\mathbf{V}} \tilde{
\mathbf{z}}_i \tilde{\mathbf{z}}_i'
\hat{\mathbf{V}}' \mathbf{E}^{kl} \hat{\mathbf{V}} \tilde{\mathbf
{z}}_i\tilde{\mathbf{z}}_i' \hat{
\mathbf{V}}' \bigr)
\]
and
\[
\tilde{\mathbf{T}}^{kl}= \operatorname{vec} \bigl(\tilde{
\mathbf{B}}^{kl} \bigr)=\tilde{\bolds{\beta}} \operatorname{vec} \bigl(
\mathbf {E}^{kl} \bigr)
\]
denote the standardized counterparts of $\hat{\mathbf{B}}^{kl}$ and
$\hat
{\mathbf{T}}^{kl}$, respectively. Then
\begin{eqnarray*}
&& \sqrt{n} \bigl(\tilde{\mathbf{T}}^{kl}-\operatorname{vec} \bigl(
\mathbf{B}^{kl} \bigr) \bigr)
\\
&&\quad= \sqrt{n} \bigl(\hat{\mathbf{T}}^{kl}- \operatorname{vec}
\bigl(\mathbf{B}^{kl} \bigr) \bigr)
\\
&&\qquad{}+ \bigl[ \bigl(\sqrt{n} (\hat{\mathbf{V}}-\mathbf {I}_p)
\otimes\mathbf{I}_p \bigr)
\\
&&\hspace*{15pt}\qquad{}+ \bigl(\mathbf{I}_p\otimes\sqrt{n} (\hat {
\mathbf{V}}-\mathbf{I}_p) \bigr) \bigr] \operatorname{vec} \bigl(
\mathbf{B}^{kl} \bigr)
\\
&&\qquad{}+\beta \bigl[ \bigl(\sqrt{n} (\hat{\mathbf{V}}- \mathbf{I}_p)
\otimes\mathbf{I}_p \bigr)
\\
&&\hspace*{21pt}\qquad{}+ \bigl(\mathbf{I}_p\otimes\sqrt{n} (\hat {
\mathbf{V}}-\mathbf{I}_p) \bigr) \bigr]\operatorname{vec} \bigl(
\mathbf{E}^{kl} \bigr).
\end{eqnarray*}
It turns out that for the asymptotics of $\hat{\mathbf{W}}$, we only need
\renewcommand{\theequation}{\arabic{equation}}
\setcounter{equation}{5}
\begin{eqnarray}
&& \sqrt{n} \bigl(\tilde{\mathbf{B}}^{kk}- \mathbf{B}^{kk}
\bigr)_{kl}
\nonumber
\\
\label{Crr} &&\quad=\sqrt{n} \bigl(\hat{\mathbf{B}}^{kk}-
\mathbf{B}^{kk} \bigr)_{kl}+ 3\sqrt{n} (\hat{\mathbf{V}}-
\mathbf{I}_p)_{kl}
\\
&&\qquad{}+(\kappa_k+3)\sqrt{n} (\hat{\mathbf{V}}-
\mathbf{I}_p)_{lk}
\nonumber
\end{eqnarray}
and
\renewcommand{\theequation}{\arabic{equation}}
\setcounter{equation}{6}
%
%
\begin{eqnarray}
&&\sqrt{n} \bigl(\tilde{\mathbf{B}}^{ll}-\mathbf{B}^{ll}
\bigr)_{kl}
\nonumber
\\
\label{Css} &&\quad=\sqrt{n} \bigl(\hat{\mathbf{B}}^{ll}-
\mathbf{B}^{ll} \bigr)_{kl}+ 3\sqrt{n} (\hat{\mathbf{V}}-
\mathbf{I}_p)_{lk}
\\
&&\qquad{}+(\kappa_l+3)\sqrt{n} (\hat{\mathbf{V}} - \mathbf
{I}_p)_{kl}.
\nonumber
\end{eqnarray}
Next, note that in the JADE procedure the matrices to be diagonalized are
\[
\tilde{\mathbf{C}}^{kl}=\tilde{\mathbf{B}}^{kl}-
\mathbf{E}^{kl}- \mathbf{E}^{lk}-\operatorname{tr} \bigl(
\mathbf{E}^{kl} \bigr) I_p,\quad k,l=1,\ldots,p.
\]
As $\sqrt{n} (\operatorname{vec}(\hat{\mathbf
{C}}^{kl})-\operatorname
{vec}(\mathbf{C}^{kl}))$
are asymptotically normal with mean zero and $\mathbf{C}^{kl}=0$, for
$k\neq l$,
and $\mathbf{C}^{kk}=\kappa_k \mathbf{E}^{kk}$, then by Lemma~\ref
{JointDiag},
$\sqrt{n} u_{kl}$ reduces to
%
%
\renewcommand{\theequation}{\arabic{equation}}
\setcounter{equation}{7}
\begin{eqnarray}
\sqrt{n} \hat{u}_{kl} &=& \frac{\kappa_k \sqrt{n} \tilde
c_{kl}^{kk}-\kappa_l
\sqrt{n} \tilde c_{kl}^{ll}}{\kappa_k^2+\kappa_l^2}+o_{P}(1)
\nonumber
\\[-8pt]
\label{uij}
\\[-8pt]
\nonumber
&=&\frac{\kappa_k \sqrt{n} \tilde b_{kl}^{kk}-
\kappa_l \sqrt{n} \tilde b_{kl}^{ll}}{\kappa_k^2+\kappa_l^2}+o_{P}(1),
\end{eqnarray}
where $\tilde c_{kl}^{kk}=(\tilde{\mathbf{C}}^{kk})_{kl}$ and $\tilde
b_{kl}^{kk}=(\tilde{\mathbf{B}}^{kk})_{kl}$.
So, asymptotically, all the information is in the matrices $\tilde
{\mathbf{B}}^{kk}$, $k=1,\ldots,p$.
As $\hat{\mathbf{W}}=\hat{\mathbf{U}}\hat{\mathbf{V}} $, where
$\hat{\mathbf{U}}$ and
$\hat
{\mathbf{V}}$ are the rotation
matrix
and the whitening matrix, respectively, we have that
\begin{eqnarray*}
\sqrt{n} (\hat{\mathbf{W}}-\mathbf{I}_p)&=& \sqrt{n} (\hat {
\mathbf{U}}\hat{\mathbf{V}}- \mathbf{I}_p)
\\
&=&\sqrt{n} (\hat{\mathbf{U}}-\mathbf{I}_p)+ \sqrt{n} ( \hat{
\mathbf{V}}-\mathbf{I}_p)+o_{P}(1).
\end{eqnarray*}
The asymptotics of the regular JADE unmixing matrix is then obtained with
$\hat{\mathbf{V}}=\hat{\mathbf{S}}^{-1/2}$, where $\hat{\mathbf
{S}}$ is the sample
covariance matrix.

Notice first that
\[
\sqrt{n} \bigl(\hat{\mathbf{S}}^{-1/2}-\mathbf{I}_p
\bigr)=-1/2 \sqrt{n} (\hat{\mathbf{S}}-\mathbf{I}_p)+o_{P}(1).
\]
Then substituting~(\ref{Crr}) and~(\ref{Css}) into~(\ref{uij}), we
have that,
for $k\neq l$,
\begin{eqnarray*}
&&\sqrt{n} \hat{w}_{kl}
\\
&&\quad=\frac{\kappa_k\sqrt{n} \hat{r}_{kl}-\kappa
_l\sqrt{n} \hat{r}_{lk}+
(3\kappa_l-3\kappa_k-\kappa_k^2)\sqrt{n} \hat{s}_{kl}}{\kappa
_k^2+\kappa_l^2}
\\
&&\qquad{}+o_{P}(1).
\end{eqnarray*}
For the diagonal elements we have simply
\[
\sqrt{n} \hat{w}_{kk}=-1/2 (\sqrt{n} \hat{s}_{kk}-1)+o_{P}(1).
\]
\upqed
\end{pf*}
\end{appendix}

\section*{Acknowledgments}
The authors wish to thank the editors and the referees for valuable
comments and suggestions.

Research supported in part by the Academy of Finland
(Grants 251965, 256291 and 268703).

%

%



\end{document}